\documentclass{article}

\usepackage{amsmath,amssymb}
\usepackage{graphs}

\newtheorem{theorem}{Theorem}[section]
\newtheorem{proposition}[theorem]{Proposition}
\newtheorem{lemma}[theorem]{Lemma}
\newtheorem{corollary}[theorem]{Corollary}

\newcommand{\countnow}{\refstepcounter{theorem}{\bf\arabic{section}.\arabic{theorem}.} }
\newcommand{\mydefinition}{\medskip{\bf Definition }\countnow}
\newcommand{\myremark}{\medskip{\bf Remark }\countnow}
\newcommand{\myexample}{\medskip{\bf Example }\countnow}

\newcommand{\N}{\mathbb N}
\newcommand{\Z}{\mathbb Z}
\newcommand{\f}{\varphi}
\newcommand{\ep}{\hfill$\square$}
\newcommand{\proc}[1]{{\bf #1}}
\begin{document}

\title{\bf Aperiodic substitutional systems and  their Bratteli diagrams}

\author{{\bf S. Bezuglyi}\\ 
Institute for Low Temperature Physics, Kharkov, Ukraine\\ 
{ bezuglyi@ilt.kharkov.ua}\\
\and{\bf J. Kwiatkowski}\\ Nicolaus Copernicus University, Toru\'{n}, Poland \\
{jkwiat@mat.uni.torun.pl}\\
\and{\bf K. Medynets\footnote{The third named author was supported by INTAS YSF-05-109-5315.}}\\
Institute for Low Temperature Physics, Kharkov, Ukraine\\
medynets@ilt.kharkov.ua}

\date{}

\maketitle




\begin{abstract}
 In the paper we study aperiodic substitutional dynamical systems arisen from non-primitive substitutions.
 We prove that the Vershik homeomorphism $\varphi$ of a stationary ordered Bratteli diagram is homeomorphic
 to an aperiodic substitutional system if and only if no restriction of $\varphi$ to a minimal component is homeomorphic to
 an odometer. We also show that every aperiodic substitutional system generated by a substitution with nesting property is homeomorphic
  to the Vershik map of a stationary ordered Bratteli diagram. It is proved that every aperiodic substitutional system is recognizable.
  The classes of  $m$-primitive substitutions and associated to them derivative substitutions are studied.
   We discuss also the notion of expansiveness for Cantor dynamical systems of finite rank.
\end{abstract}

\section{Introduction}
During last decade, minimal homeomorphisms of a Cantor set (Cantor
minimal systems, in other words)
have been thoroughly studied in many papers. The most powerful tool in the study of Cantor minimal systems
is the concept of Bratteli diagrams. It was shown in \cite{herman_putnam_skau:1992}
that every minimal homeomorphism is topologically conjugate to the Vershik map of an ordered simple Bratteli diagram.
Such a realization of minimal homeomorphisms allowed one to prove many deep results clarifying properties and orbit structure of minimal Cantor systems.
   We mention here the study of orbit equivalence and affability of Cantor minimal systems which was conducted in papers
    by Glasner, Giordano, Matui, Putnam, Skau, Weiss, and others (see \cite{giordano_putnam_skau:1995}, \cite{giordano_putnam_skau:1999}, \cite{giordano_putnam_skau:2004},
    \cite{giordano_matui_putnam_skau:2006}, \cite{glasner_weiss:1995}).
 Another meaningful usage of Bratteli diagrams was given in the papers \cite{forrest:1997} and \cite{durand_host_scau:1999}.
  They answered the natural question of the description of the class of minimal homeomorphisms which can be represented by Bratteli diagrams
   of the simplest form, i.e., by stationary Bratteli diagrams. It turns out that this class is constituted by minimal substitutional systems and odometers.

Motivated by these remarkable achievements, we are interested in
the following question: Is the assumption of minimality crucial in
proving
 these results?  In the paper \cite{bezuglyi_dooley_medynets:2005}, we considered aperiodic Cantor systems and proved
 the existence of Kakutani-Rokhlin partitions for them. In contrast to the minimal case,
 we cannot start with an arbitrary clopen set to produce a Kakutani-Rokhlin partition.
Nevertheless, it was proved in \cite{medynets:2006} that, given an
aperiodic homeomorphism $T$, there exists a sequence of nested
Kakutani-Rokhlin partitions which generate the topology. Thereby
we showed that every aperiodic homeomorphism can be realized as
the Vershik map of an ordered Bratteli diagram      (see details
in Section 2). However, the structure of Bratteli diagrams for
aperiodic homeomorphisms is still unclear --- in comparison with
the minimal case where each simple Bratteli diagram can be given
an order that defines a continuous Vershik map
\cite{herman_putnam_skau:1992},  not every Bratteli
diagram admits, in general, such an ordering \cite{medynets:2006}.

The primary goal of this work is the study of aperiodic
(non-minimal) substitutional systems and finding  explicit descriptions of their Bratteli-Vershik models.  To the best our knowledge, non-primitive substitutions have not been systematically studied yet. In our study we are mostly based on two articles.
 The first one is the paper by Durand, Host, and Skau  \cite{durand_host_scau:1999} where the Bratteli diagrams for primitive substitutional dynamical systems were thoroughly studied. The second one is the very recent work by Downarowicz and Maass \cite{downarowich_maass:2006} where the authors  suggested a very
fruitful idea of coding of dynamics by means of the so-called
$j$-symbols. In other words, this approach gives  a symbolic
interpretation of the technique of Kakutani-Rokhlin partitions and
Bratteli diagrams. The basic references to the study of substitutional dynamical systems are \cite{fogg:book}
and \cite{queffelec:book}. We also refer the reader to the book
\cite{kurka:book} for a comprehensive exposition of the symbolic
dynamics.

We will use the standard notation of the theory of substitutional systems. Denote by $A$  a finite alphabet and by $A^+$ the set of all non-empty words over $A$. Let $\sigma :A\rightarrow A^+$ be a substitution. We denote by $A_l$ the set of all letters $a\in A$ such that $|\sigma^n(a)|\to\infty$ as $n\to\infty$. Let $A_s=A\setminus
A_l$. We say that a substitution $\sigma$  has a {\it nesting
property} if at least one of the following two conditions holds: (1) for every $a\in A_l$ the word $\sigma(a)$ starts with a letter from $A_l$; (2) for every $a\in A_l$ the
word $\sigma(a)$ ends with a letter from $A_l$. Clearly, if $A =
A_l$ then $\sigma$ has the nesting property.

\medskip\noindent The paper is organized as follows:

\medskip\noindent {\it Section \ref{Section_NotionBratteliDiagram}: Bratteli-Vershik models of Cantor
aperiodic systems.} In the section, we consider aperiodic
homeomorphisms of a Cantor set and discuss the notions of an
ordered  Bratteli diagram and the Vershik map associated to
aperiodic Cantor systems. We  outline the proof of the fact that
any aperiodic homeomorphism of a Cantor set is conjugate to the
Vershik map of an ordered Bratteli diagram (see
\cite{medynets:2006}). This result is a foundation stone for our
further  research.

\medskip\noindent {\it Section \ref{SectionSymbolicFormalism}: Downarowicz-Maass' symbolic representation.}
 In the paper  \cite{downarowich_maass:2006} the
authors suggested a method of coding of  dynamics on Bratteli
diagrams by means of {\it $j$-symbols}. In the section, we
generalize ideas and results from \cite{downarowich_maass:2006}
and give  abstract definitions of $j$-symbols, $j$-sequences,
etc. They can be used  to study  dynamics of different nature, for instance, Bratteli-Vershik systems and substitutional dynamical systems. The
main advantage of this approach is that it allows one to use the
machinery of symbolic dynamics for solving some problems of Cantor
dynamics.
 The technique used in this section (see Propositions \ref{PropositionSymbolicFormasimNoCommonCuts} and \ref{PropositionInfectionLemma})
 is applied in the proofs of our main results.

\medskip\noindent{\it Section \ref{SectionExpansiveFiniteRank}: Finite rank aperiodic Cantor systems}.
 We apply in this section the technique of $j$-sequences
developed in Section \ref{SectionSymbolicFormalism} to the study
of Cantor aperiodic systems of  finite rank. Generalizing the main
result of \cite{downarowich_maass:2006}, we show that the  Vershik
map of an ordered Bratteli diagram with uniformly bounded number
of vertices at each level whose restrictions to minimal components
are not conjugate to odometers is expansive, i.e. this map is
homeomorphic to a subshift over a finite alphabet.

\medskip\noindent {\it Section \ref{SectionNotionOfSubstDynamSystem}:
Recognizability of aperiodic substitutional systems.} We discuss
in this section the properties of an arbitrary substitution $\sigma$
over a finite alphabet $A$ and the substitutional dynamical system
$(X_\sigma,T_\sigma)$ defined by $\sigma$ (the rigorous definitions are given  in Section  \ref{SectionNotionOfSubstDynamSystem}). We prove that the
number of minimal components of the system $(X_\sigma,T_\sigma)$
cannot exceed $|A|$. We also show that for an arbitrary
substitution $\sigma$ any point $x\in X_\sigma$ can be written as
\begin{equation}\label{EquationRealizationX}
 x=T^i_\sigma\sigma(y)\mbox{ for some }y\in X_\sigma\mbox{ and some
}i=0,\ldots,|\sigma(y[0])|-1.
 \end{equation}
If $\sigma$ is a primitive aperiodic substitution,
then representation (\ref{EquationRealizationX}) is unique for
each $x$. This fact is a consequence of recognizability property of
primitive substitutions established by Moss\'e
\cite{mosse:1992,mosse:1996}. We generalize this result by showing
that if we omit the condition of primitivity for $\sigma$, then
the uniqueness  of the representation (\ref{EquationRealizationX}) still holds.  The proof of
the result involves the ideas of  Downarowicz-Maass' symbolic
representation \cite{downarowich_maass:2006} and the technique
developed in Section  \ref{SectionSymbolicFormalism}.

\medskip\noindent {\it Section
\ref{SectionDiagramsForSubstitutions}: Stationary Bratteli-Vershik
models  vs. aperiodic substitutional systems.} In the section we
show that every expansive Vershik map of a stationary ordered
Bratteli diagram is homeomorphic to a substitutional dynamical
systems. Conversely, we show that every aperiodic substitutional
dynamical system constructed by a substitution satisfying the {\it nested} property is homeomorphic   to the Vershik map of a
stationary ordered Bratteli diagram. These results generalize
those proved in \cite{durand_host_scau:1999} for minimal
(primitive) substitutional systems.

\medskip\noindent {\it Section \ref{Section_DerivativeSubstitutions}: Derivative substitutions.}
The section is devoted to the study of derivative substitutions
associated with aperiodic substitutions. This notion was first considered in
\cite{durand_host_scau:1999} and \cite{durand:1998} for primitive
substitutions. We show how using derivative substitutions one can
find a Bratteli-Vershik realization of a substitutional system.
This approach differs from that of Section
\ref{SectionDiagramsForSubstitutions} and sometimes gives a
`simpler' Bratteli diagram.

\medskip\noindent{\it Appendix: Description of the phase space
$X_\sigma$.} In the appendix we give an explicit description of
elements from the space $X_\sigma$ and establish some
supplementary results on general substitutional  systems.

%
%

\section{Bratteli-Vershik Models of Cantor Aperiodic
Systems}\label{Section_NotionBratteliDiagram}

In this section we show how one can associate a Bratteli-Vershik
dynamical system to every Cantor aperiodic (non-minimal)
homeomorphism. This result was announced  in
\cite{medynets:2006}. Since the notion of Bratteli diagrams has
been discussed in many well known  papers on Cantor dynamics (e.g.
\cite{herman_putnam_skau:1992}, \cite{giordano_putnam_skau:1995},
\cite{giordano_putnam_skau:2004}), we give only the main steps of
the construction.

By a {\it Cantor set} $X$ we understand any zero-dimensional
compact metric space without isolated points. Recall that a
homeomorphism $T$ is called {\it aperiodic} if
every $T$-orbit is infinite; a homeomorphism $T: X\rightarrow X$ is called {\it minimal} if every orbit of $T$ is dense in $X$.

\mydefinition \label{Delete_1} A {\it Bratteli diagram} is an
infinite graph $B=(V,E)$ such that the vertex set
$V=\bigcup_{i\geq 0}V_i$ and the edge set $E=\bigcup_{i\geq 1}E_i$
are
partitioned into disjoint sets $V_i$ and $E_i$ such that\\
(i) $V_0=\{v_0\}$ is a single point; \\
(ii) $V_i$ and $E_i$ are finite sets;
\\
(iii) there exist a range map $r$ and a source map $s$ from $E$ to
$V$ such that $r(E_i)= V_i$, $s(E_i)= V_{i-1}$, and $s^{-1}(v)\neq\emptyset$, $r^{-1}(v')\neq\emptyset$ for all $v\in V$ and $v'\in V\setminus V_0$.
 \medbreak

The pair  $(V_i,E_i)$ is called the $i$-th level of the diagram $B$. We write
$e(v,v')$ to denote an edge $e$ such that $s(e)=v$ and $r(e)=v'$.

A finite or infinite sequence of edges $(e_i : e_i\in E_i)$ such
that $r(e_{i-1})=s(e_i)$ is called a finite or infinite path,
respectively. It follows from the definition that every vertex
$v\in V $ is connected to $v_0$ by a finite path and the set
$E(v_0,v)$ of all such paths  is finite. For a Bratteli diagram
$B$, we denote by $X_B$ the set of infinite paths. We endow
the set $X_B$ with the topology generated by  cylinder sets
$U(e_1,\ldots,e_n)=\{x\in X_B : x_i=e_i,\;i=1,\ldots,n\}$, where
$(e_1,\ldots,e_n)$ is a finite path of $B$. The set $X_B$ is a 0-dimensional compact metric space with
respect to this topology.

\myremark In general, the set $X_B$ may have isolated
points. We  do not assume that $X_B$ is a perfect space; our results remain true for any 0-dimensional compact metric space.
\medbreak

Let $B=(V,E,\leq)$ be a Bratteli
diagram $(V,E)$ equipped with a partial order $\leq$ defined on
each $E_i,\ i=1,2,...,$ such that edges $e,e'$ are comparable if
and only if $r(e)=r(e')$; in other words, a linear order $\leq$ is
defined on each (finite) set $r^{-1}(v),\ v\in V\setminus V_0$.
For a Bratteli diagram $(V,E)$ equipped with such a partial order
$\leq$ on $E$, one can also define a partial lexicographic order
on the set $E_{k+1}\circ\cdots\circ E_l$ of all paths from $V_k$
to $V_l$: $(e_{k+1},...,e_l) < (f_{k+1},...,f_l)$ if and only if
for some $i$ with $k+1\le i\le l$, $e_j=f_j$ for $i<j\le l$ and
$e_i< f_i$. Then any two paths from $E(v_0,v)$ are comparable with
respect to the introduced lexicographic order. We call a path $e=
(e_1,e_2,..., e_i,...)$ {\it maximal (minimal)} if every $e_i$ is
maximal (minimal) amongst all elements from $r^{-1}(r(e_i))$. Notice that there are unique minimal and
maximal paths in $E(v_0,v)$ for each $v\in V_i,\ i\ge 0$.

Denote the sets of all maximal and minimal paths in $X_B$ by
$X_{max}$ and $X_{min}$, respectively. It is not hard to see that
$X_{max}$ and $X_{min}$ are  non-empty closed sets.

\mydefinition A Bratteli diagram $B=(V,E)$ together with a partial
order $\leq$ on $E$ is called {\it an ordered Bratteli diagram
$B=(V,E,\leq)$}. \medbreak

\mydefinition\label{Defintion_VershikMap}  Let $B=(V,E,\leq)$ be
an ordered Bratteli diagram.  We say that a homeomorphism
$\varphi_B : X_B\rightarrow X_B$ is a {\it Vershik homeomorphism (map)}
if it satisfies the following conditions:
\\
(i) $\varphi_B(X_{max})=X_{min}$;
 \\
(ii) if $x=(x_1,x_2,\ldots)\notin X_{max}$, then
$$\varphi_B(x_1,x_2,\ldots)=(x_1^0,\ldots,x_{k-1}^0,\overline
{x_k},x_{k+1},x_{k+2},\ldots),$$ where $k=\min\{n\geq 1 :
x_n\mbox{ is not maximal}\}$, $\overline{x_k}$ is the successor of
$x_k$ in the set $r^{-1}(r(x_k))$, and $(x_1^0,\ldots,x_{k-1}^0)$
is the minimal path in $E(v_0,s(\overline{x_k}))$. \medbreak

It is well-known that every simple  Bratteli diagram $B$ has an ordering such that the sets $X_{\max}$ and $X_{\min}$ are
singletons. In contrast to this case,  non-simple Bratteli diagrams may not have such an ordering, in general. We notice also that not every order of a Bratteli diagram defines a continuous Vershik map $\varphi_B$. Moreover, it was
shown in \cite{medynets:2006} that {\it there exists a Bratteli
diagram $B$ such that any ordering of $B$ does not produce the
(continuous) Vershik map.}

Let $B=(V,E,\leq)$ be an ordered Bratteli diagram such that the
Vershik map $\varphi_B$ exists. Then the pair $(X_B,\varphi_B)$ is
called a {\it Bratteli-Vershik system}.

To construct a Bratteli diagram associated to a homeomorphism of a
Cantor set, one needs to work with sequences of Kakutani-Rokhlin
partitions. For a  minimal homeomorphism,  such sequences can be
easily produced via the first return functions for any clopen
sets.

We consider now aperiodic Cantor systems and show how we can
construct Bratteli diagrams in this case
\cite{bezuglyi_dooley_medynets:2005}, \cite{medynets:2006}.

\mydefinition Let $(X,T)$ be a Cantor dynamical system. (1) By a
{\it Kakutani-Rokhlin (K-R) partition} $\mathcal P$ of $X$, we
understand any partition of $X$ into clopen sets of the form
$$
\mathcal P=\{T^iA(v)\;:\;v=1,\ldots,m; \ \ i=0,\ldots,h(v)-1\}.
$$
The {\it base} of the partition  $\mathcal
P$ is $B(\mathcal P)=\bigcup_{v=1}^m A(v)$. A family of disjoint
sets $\xi (v)= \{A(v),TA(v),\ldots,T^{n-1}A(v)\}$ is called a {\it
$T$-tower} of {\it height} $h(\xi)= n$. Then $\mathcal P$ is the disjoint union of clopen $T$-towers $\xi(v)$.

(2) A sequence of K-R partitions $\{\mathcal P_n\}$ is called {\it
nested} if $\mathcal P_{n+1}$ refines $\mathcal P_n$ and
$B(\mathcal P_{n+1})\subseteq B(\mathcal P_n)$ for any $n$.
\medbreak

For a Cantor aperiodic system $(X, T)$, a clopen set  $A$ is
called a {\it complete $T$-section} if $A$ meets every $T$-orbit
at least once. A point $x \in A$ is called {\it recurrent} with
respect to $A$ if there exists $n \in \N$ such that $T^nx \in A$.
By compactness of $X$, every clopen complete
$T$-section $A$ consists of recurrent points. Thus, using the
first return function $n_A(x) = \min \{n \in  \N : T^nx \in A\}$
for a complete $T$-section $A$, we can construct a K-R partition
$\mathcal P$ of $X$ with $B(\mathcal P) =A$.

\begin{theorem} \cite{medynets:2006} \label{Theorem-ExistenceOF_K-R_Partitions} Let $(X, T)$ be a Cantor aperiodic system. There exists a sequence of K-R partitions $\{\mathcal P_n\}$ of $X$ such that for all $n\ge 1$: (i) $\mathcal P_{n+1}$ refines $\mathcal P_n$; (ii) $h_{n+1} > h_n$, where $h_n = \min \{h(\xi) : \xi \in  \mathcal P_n$\}; (iii) $B(\mathcal P_{n+1})\subseteq B(\mathcal P_n)$; (iv) the elements of  partitions $\{\mathcal P_n\}$ generate the clopen topology on $X$.
\end{theorem}

Let $(X, T)$ be a Cantor aperiodic system. We say that a closed
set $Y \subset X$  is a {\it basic} set if every clopen
neighborhood of $Y$ is a complete $T$-section and $Y$ meets every
$T$-orbit at most once. If a sequence of K-R partitions
$\{\mathcal P_n\}$ satisfies the conditions of Theorem
\ref{Theorem-ExistenceOF_K-R_Partitions}, then $Y = \bigcap_n
B(\mathcal P_n)$ is a basic set.

\begin{theorem} \cite{medynets:2006} \label{ExistenveBratteliDiagrams}
 Let $(X, T)$ be a Cantor aperiodic system with a basic set $Y$. There exists an ordered Bratteli diagram
 $B = (V,E,\leq)$ such that $(X, T)$ is conjugate to a Bratteli-Vershik model $(X_B, \varphi_B)$. The homeomorphism implementing the conjugacy between $T$ and
$\varphi_B$ maps the basic set $Y$ onto the set $X_{\min}$ of all
minimal paths of $X_B$.
\end{theorem}

Our primary goal is to construct Bratteli-Vershik models for
aperiodic substitutional systems. We recall briefly the main steps of
this construction for an aperiodic homeomorphisms of a Cantor set.

Let  $(X,T)$ be a Cantor aperiodic system. By Theorem
\ref{Theorem-ExistenceOF_K-R_Partitions} find a sequence of nested
K-R partitions: $\mathcal P_0 = X$,
\begin{equation}\mathcal P_n=\{T^iA(v,n)\; :\; v=1,\ldots,m(n); \ \  i=0,\ldots,h(v,n)-1\},\;n\geq 1,
\end{equation}
that generates the clopen topology on $X$. Set
$$\xi(v,n)=\{A(v,n),\ldots,T^{h(v,n)-1}A(v,n)\}\mbox{ for
}v=1,\ldots,m(n).
$$
Define an ordered Bratteli diagram $B=(V,E,\leq)$ as follows:

(i) Define the set of vertices by $V_0=\{v_0\}$ and
$V_n=\{1,\ldots,m(n)\}$ for  $n\geq 1$.

(ii) Define the set  of edges $E_n$  between the consecutive
levels $V_{n-1}$ and $V_n$ by the incidence matrix
$M(n)=\{m_{vw}(n) : v\in V_{n},\;w\in V_{n-1}\}$, where
$$
m_{vw}(n)= |\{0\leq i<h(v,n) : T^iA(v,n)\subset A(w,n-1)\}|.
$$
In other words, we fix a vertex $v\in V_n$ and define $V(v,n-1)$
as the set of all vertices from $V_{n-1}$ such that $\xi(v,n)$
intersects $\xi(w,n-1)$, and a vertex $w$ appears in $V(v,n-1)$ as
many times as $\xi(v,n)$ intersects $\xi(w,n-1)$. Then we connect
$v$ to each vertex $w\in V(v,n-1)$ taking into account the
multiplicity of appearance of $w$ in $V(v,n-1)$.

(iii) To define the ordering $\leq$ on $E$ we take the clopen set
$A(v,n)$. Then tracing the orbit of $A(v,n)$ within the $T$-tower
$\xi(v,n)$, we see that $T^iA(v,n)$ consecutively meets the sets
$A(w_1,n-1)$, $A(w_2,n-1)$,..., $A(w_{k_v},n-1)$ (some of them can
occur several times). This defines the set of edges $r^{-1}(v)$.
Enumerate the edges from $r^{-1}(v)$ as follows:
$e(w_1,v)<e(w_2,v)<\ldots<e(w_{k_v},v)$.

Since the partitions $\{\mathcal P_n\}$ generate the topology of
$X$, for each point $x\in X$ there is a unique sequence
$i(x)=\{(v_n,i_n)\}_{v\in V_n;\,0\leq i_n<h(v_n,n)}$ such that
$$
\{x\}=\bigcap_{n=1}^\infty T^{i_n}A(v_n,n).
$$
Define the map $\theta : X\rightarrow X_B$ by
$$\theta(x)=\bigcap_{n\geq 1}U(y_1,\ldots,y_n)
$$ where
$(y_1,\ldots,y_n)$ is the $i_n$-th finite path in $E(v_0,v_n)$
with respect to the lexicographical ordering on $E(v_0,v_n)$.  It can be easily checked that $\theta$ is a homeomorphism.

Define $\varphi_B=\theta \circ T\circ\theta^{-1}$. Clearly, the homeomorphism  $\varphi_B$ satisfies Definition
\ref{Defintion_VershikMap} with the ordering $\leq$. Thus, we
obtain that $(X,T)$ is conjugate  to the Bratteli-Vershik system
$(X_B,\varphi_B)$.

\medskip {\bf Open question.} It is well-known that for  every simple properly ordered Bratteli diagram $B=(V,E)$ there exists the continuous  Vershik map, see \cite{herman_putnam_skau:1992}. On the other hand there are Bratteli diagrams such that any ordering on them will not lead to a {\it continuous} Vershik map \cite{medynets:2006}.
The question is: How can one  describe the class of Bratteli
diagrams which admit continuous  Vershik maps on the set of
infinite paths?

%
%
%

\section{Downarowicz-Maass' Symbolic Representation}\label{SectionSymbolicFormalism}

In the section, we develop the ideas of the paper
\cite{downarowich_maass:2006} and apply them below to solving the
following problems: recognizability of aperiodic substitutions
(Section \ref{SectionNotionOfSubstDynamSystem}) and expansiveness
of Vershik maps (Section \ref{SectionExpansiveFiniteRank}). The
results stated in this section are proved in
\cite{downarowich_maass:2006}.  Since they are of crucial
importance for us, we reproduce their proofs here with minor alterations.

Fix a sequence of finite alphabets $\{A_i\}_{i\geq 0}$.

\mydefinition (1) Take $a\in A_j$, $j\geq 1$. By a {\it
$j$-symbol} $[a]_j$ we understand a finite matrix with $j+1$ rows
numbered from $\{0,1,\ldots,j\}$ which has the following
structure:

(i) The row $j$ (the bottom row of $[a]_j$) consists of one box
carrying the symbol $a$ (in other words, labeled by $a$) whose
length is extending over the full width of the matrix.

(ii) The row $j-1$ consists of $m_{j-1}$ boxes labeled by letters
$b_{0},\ldots,b_{m_{j-1}-1}$ from the alphabet $A_{j-1}$. The
total length of these boxes is equal to the width of the matrix
$[a]_j$.

(iii) The row $j-2$ consists of $m_{j-1}$ group of boxes. The
$i$-th group, $i= 0,\ldots,m_{j-1}-1$, is located exactly over the
box of $(j-1)$-th row labeled by $b_i$; the total length of boxes
from the group $i$ is equal to the length of the box $b_i$. Each
box from the row $j-2$ is labeled by a letter from the alphabet
$A_{j-2}$.

(iv) All rows above have the same structure. The row 0 (the first
row of the matrix) consists of $m_0$ boxes of length 1 which are
labeled by letters from $A_0$.

(2) A {\it 0-symbol} is  always a $1\times 1$ matrix viewed as a
box of length 1 which is labeled by a letter from $A_0$.
\medbreak

We observe  that the set of all
$j$-symbols, $j\geq 1$, is infinite. Denote by $\mathcal A_j$ any
finite set of $j$-symbols over the alphabets $A_0,\cdots, A_j$, $j\geq 0$. When we use the notation $\mathcal A_j$, it will be clear from the
context which set of $j$-symbols we mean.

\mydefinition\label{Definition_Agreeable_J-Symbols} Let  maps
$\sigma_i : A_i\rightarrow A_{i-1}^+$, $i\geq 1$, be given. For
every $i$ define by induction the family of $i$-symbols
$\{\mathcal A_i\}_{i\geq 0}$:

(1) $\mathcal A_0 = \{[a]_0 : a\in A_0\}$.

(2) Suppose that the family $\mathcal A_{i-1}$ is defined.  For
$a\in A_i$, let $\sigma_i(a) = a_0\cdots a_m$ where
$a_0,\ldots,a_m \in A_{i-1}$. Define the $i$-symbol $[a]_i$ as the
matrix whose first $i$ rows are concatenations of $(i-1)$-symbols
$[a_0]_{i-1},\ldots,[a_{m}]_{i-1}$ from $\mathcal A_{i-1}$. Then
the family $\mathcal A_i$ is formed by $i$-symbols $\{[a]_i : a\in
A_i\}$.

We call this families of  $i$-symbols $\{\mathcal A_i\}_{i\geq 0}$ {\it agreeable}.\medbreak

\myexample \label{Example_2_symbol} Let $A_2=\{a,b,...\}$,
$A_1=\{x,y...\}$, $A_0=\{\alpha, \beta,...\}$. Then the following
matrix gives an example of a 2-symbol $[a]_2$: \medbreak

\begin{center}
{\small \begin{tabular}{ | c | c | c | c |} \hline  $\alpha$ & $\beta$ & $\gamma$ & $\delta$\\
 \hline \multicolumn{2}{|c|} x & \multicolumn{2}{|c|} y \\ \hline \multicolumn{4}{|c|}
 a \\ \hline
  \end{tabular}}
\end{center}

In this example, we have that $\sigma_2(a) = xy$, $\sigma_1(x)
=\alpha\beta$, and $\sigma_1(y) = \gamma\delta$ where $\sigma_i:
A_i \to A^+_{i-1}$, $i= 1,2$. \medbreak

\myremark Suppose that the families of $j$-symbols $\{\mathcal A_j\}_{j\geq
0}$ are agreeable. Then any row $i\leq j$ of a $j$-symbol
$[a]_j\in \mathcal A_j$ completely determines all rows above.
Furthermore, it follows from Definition
\ref{Definition_Agreeable_J-Symbols} that each letter $a\in A_j$
is uniquely assigned to the $j$-symbol $[a]_j$ from $\mathcal
A_j$. \medbreak

We give now two principal constructions of agreeable families.

\myexample \label{j-symbols in Bratteli diagrams} Let
$B=(V,E,\leq)$ be an ordered Bratteli diagram with the path space
$X_B$. Set $A_i= V_i$, $i\geq 0$. Define the map $\sigma_i :
A_{i}\rightarrow A_{i-1}^+\ (i \ge 1)$ as follows: take a vertex
$v\in V_i$, write down all edges $e_1<e_2<\ldots<e_m$ from
$r^{-1}(v)$ with respect to the ordering $\leq$, and define
$$
\sigma_i(v)=v_1\ldots v_m
$$
where $v_k = s(e_k)$, $k=1,\ldots,m$. Then the alphabets
$\{A_i\}_{i\geq 0}$ and maps $\sigma_i : A_i\rightarrow A_{i-1}^+$
define the sets $\{\mathcal A_j\}_{j\geq 0}$ of agreeable
$j$-symbols,  see Definition \ref{Definition_Agreeable_J-Symbols}.
\medbreak

To illustrate this example consider the following Bratteli diagram
with the edges ordered from the left to right.
 \unitlength=0.5mm

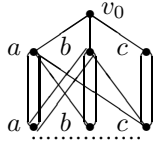
\begin{center}
\begin{picture}(40,40)
\put(20,35){\circle*{2}}
\put(23,35){$v_0$}
\put(-2,24){{\bf $a$}}
\put(5,25){\circle*{2}}
\put(12,24){{\bf $b$}}
\put(20,25){\circle*{2}}
\put(27,24){{\bf $c$}}
\put(35,25){\circle*{2}}
    \put(5,25){\line(3,2){15}}
    \put(20,25){\line(0,1){10}}
    \put(35,25){\line(-3,2){15}}
\put(-2,4){{\bf $a$}}
\put(5,5){\circle*{2}}
\put(12,4){{\bf $b$}}
\put(20,5){\circle*{2}}
\put(27,4){{\bf $c$}}
\put(35,5){\circle*{2}}
%
\put(5,15){\oval(3,20)}
\put(20,15){\oval(3,20)}
\put(35,15){\oval(3,20)}
\put(4,6){\line(3,4){15}}
\put(6,4){\line(3,4){15}}
\put(20,5){\line(-3,4){15}}
\put(35,5){\line(-3,2){30}}
\put(35,5){\line(-3,4){15}}
\multiput(5,2)(2,0){15}{\circle*{0.7}}
\end{picture}
\end{center}
Then the 2-symbol $[c]_2$ can be viewed as the matrix
\begin{center}
 {\small $$\begin{array}{ | c | c | c | c|} \hline
v_0 & v_0 & v_0 & v_0 \\ \hline a & b & c & c \\ \hline
\multicolumn{4}{|c|} c\\ \hline
 \end{array}$$}
\end{center}
\medbreak
\medskip

\myexample \label{j-symbols in substitutions} Let $\tau :
A\rightarrow A^+$ be a substitution defined on a finite alphabet
$A$. Set $A_0 = A$ and $A_j = \{\tau^j(a) : a\in A\}, j>0$. Define
$\sigma_j : A_j \to A^+_{j-1}$ by setting $\sigma_j(\tau^j(a)) =
\tau^{j-1}(a_1)\cdots \tau^{j-1}(a_m)$ where $a_1\cdots a_m =
\tau(a)$. It follows from Definition
\ref{Definition_Agreeable_J-Symbols} that the data $\{(A_j,
\sigma_j)\}_{j\geq 0}$ define the agreeable families of
$j$-symbols $\{\mathcal A_j\}_{j\geq 0}$. \medbreak

\mydefinition\label{Definition_J-sequences} Let $\{\mathcal
A_0,\ldots,\mathcal A_j\}$ be agreeable families of symbols.
Denote by $\mathcal  Z_j$ the set of all matrices with $j+1$
two-sided infinite rows which are obtained by  concatenation of
$j$-symbols. Every element of  $\mathcal Z_j$ is a sequence
$\{[a_n]_j\}_{n\in\mathbb Z}$ with $[a_n]_j\in\mathcal A_j$.
Elements of $\mathcal Z_j$ are called {\it $j$-sequences. }
\medbreak

For $x\in\mathcal Z_j$ and  $n\in\mathbb Z$, denote by $x(n)$ the
$n$-th column of $x$. The $i$-th row of the matrix $x$ ($i < j$)
is a concatenation of boxes of variable lengths which are  labeled
by $i$-symbols from $\mathcal A_i$. So it is natural to represent
$x(n)$  as a column  whose  $i$-th entry carries the following
information: the letter $a\in A_i$ such that the $i$-symbol
$[a]_i$ intersects the column $x(n)$ in $x$ and an integer $k>0$
that denotes the coordinate of the column $x(n)$ within the
$i$-symbol $[a]_i$.

This observation allows us to write down the matrix ($j$-sequence)
$x$ as the sequence $x=\{x(n)\}_{n\in\mathbb Z}$. Define $T:
\mathcal Z_j\rightarrow \mathcal Z_j$ by $(Tx)(n)=x(n + 1)$ for
all $n\in \mathbb Z$. Then $(\mathcal Z_j,T)$ is a subshift over a
finite alphabet. \medbreak

The following picture gives an example of a 2-sequence.

{\small $$\begin{array}{c|c|c|c|c|c|c|c|c|c|c|c|c|c|c} \hline
 \ldots &   x_{-6} & x_{-5} & x_{-4} & x_{-3}
 & x_{-2} & x_{-1} & x_{0} & x_{1} & x_{2} & x_{3} & x_{4}
 & x_{5} & x_{6}  & \ldots \\
 \hline  \ldots & \multicolumn{3}{|c|}  {y_{-2}} &
 \multicolumn{2}{|c|} {y_{-1}} &  \multicolumn{4}{|c|} {y_{0}} &  \multicolumn{4}{|c|}
 {y_{1}} & \ldots \\
 \hline \ldots &  \multicolumn{5}{|c|} {z_{-1}} &  \multicolumn{8}{|c|} {z_{0}}
 & \ldots \\ \hline
\end{array}$$}

\mydefinition\label{Definition_Compatible_J-sequences} Let
$\{\mathcal A_0,\mathcal A_1,\ldots\}$ be  agreeable families of symbols and
$x,y\in\mathcal Z_j$ for some $j>0$.

(1)  $x$ and $y$ are called {\it $i$-compatible}, $i\leq j$, if
the  $i$-th rows of $x$ and $y$ coincide. Observe that since the
families $\{A_k\}_{k\geq   0}$ are agreeable, the $i'$-th rows of
$x$ and $y$ also coincide for all $i'<i$.

(2) If $x\neq y$ and they are compatible, then the maximal $i\leq j$ such that the $i$-th rows of $x$ and $y$ coincide is called the {\it depth} of $x$ and $y$.

(3) If the $i$-th rows  of $x$ and $y$ are different, then $x$ and $y$
are called {\it $i$-separated}.

(4) $x$ and $y$ are said to have a {\it common $j$-cut} if there
exist  $n\in \N$ and $j$-symbols $[a]_j$ and $[b]_j$ such that
$[a]_j$ and $[b]_j$ appear in $x$ and $y$, respectively, at the
position $n$. \medbreak

\myremark (1) We note that if a pair $(x,y)$ has a common $j$-cut
at a coordinate $n$, then it has common $i$-cuts for all $i\leq j$
at the same coordinate $n$.

(2) For a $j$-symbol $[a]_j$, let $|[a]_j|$ stand for the width of
the matrix $[a]_j$. We observe that if a $j$-symbol $[a]_j$ appears
in $x\in\mathcal Z_j$ at the position $n$, then there is a
$j$-symbol $[a']_j$ that appears in $x$ at the position
$n+|[a]_j|$.

(3) Denote by $\pi_i : \mathcal Z_j\rightarrow \mathcal Z_i$,
$i<j$, the map that restricts each  $x\in \mathcal Z_j$ to its
first $(i+1)$ rows. Then $\pi_i$ is a factor map from $(\mathcal
Z_j,T)$ onto $(\mathcal Z_i,T)$.\medbreak

It is not hard to see that the technique of $j$-symbols is an
interpretation of the well-known technique of Kakutani-Rokhlin
partitions. However, the usage of $j$-symbols sometimes is more
convenient as it allows us to manipulate with blocks, symbols, and
other symbolic objects.

\begin{proposition}\cite{downarowich_maass:2006}
\label{PropositionSymbolicFormasimNoCommonCuts} Let $\{\mathcal
A_j\}_{j\geq 0}$ be  agreeable families of $j$-symbols. Suppose
also that $|\mathcal A_j|\leq K<\infty$ for all $j\geq 0$. Then
for any $n\in\mathbb N$ there exist $i,j\in\mathbb N$ with $n\leq
i<j$ such that any pair of $j$-sequences with  depth $i$ has no
common $j$-cuts.
\end{proposition}
\proc{Proof.} Assume the converse. Then take $n_0\in\mathbb N$
such that for any $i=n_0,\dots,j-1$ there exists a pair
$(x_i,y_i)$ of $j$-sequences with depth $i$ that has a common
$j$-cut. We  set $j=n_0+K$.

(I) Consider the pair $(x_{j-1},y_{j-1})$. By assumption, this
pair is $(j-1)$-compatible, $j$-separated and has a common
$j$-cut.

 Fix any common $j$-cut for $x_{j-1}$ and $y_{j-1}$, i.e. we
fix a position $n$ at which   possibly different $j$-symbols appear
in $x_{j-1}$ and $y_{j-1}$. Consider the $j$-symbols occurring at
the position $n$ in $x_{j-1}$ and $y_{j-1}$. If these $j$-symbols
are the same, then the following $j$-cut is also common and we
consider the following $j$-symbols (to the right). Let $u$ and $v$
be the first different $j$-symbols to the right of $n$. If there
are no such symbols, we analyze $x_{j-1}$ and $y_{j-1}$ to the
left of $n$.

If the $j$-symbols $u$ and $v$ have the same length, then the
following $j$-cut is also common. In this case, we do not change
anything and continue checking the following pair of $j$-symbols.

If, say, $u$ is longer than $v$, we modify the set of $j$-symbols.
We replace the $j$-symbol $u$ from $\mathcal A_j$ by  the
concatenation of two  $j$-symbols, $v\in \mathcal A_j$ and a new
symbol $u'$ which is defined as follows. The matrix $[u']_j$ is
formed by $|u| - |v|$ the right most columns of the matrix $[u]_j$
with the $(j+1)$-th row labeled by $u'$. In other words, we
substitute the last row in the symbol $u$ by two boxes: $v$ and a
new symbol $u'$ which has the complementary length. Observe that
the rows above the last one are not changed (recall that $x_{j-1}$
and $y_{j-1}$ are $(j-1)$-compatible).

Next we replace  occurrence of the
$j$-symbol $u$ in every element of $\mathcal Z_{j}$ by the
concatenation $vu'$. This procedure lead to a symbolic system
which is topologically  conjugate to $(\mathcal Z_{j}, T)$. Note
also that we produce more $j$-cuts and never remove them. If there
is still a non-common cut, we repeat the described construction
and substitute some $j$-symbol with an existing one and a new one.

As soon as we get that all $j$-cuts to the right of $n$ are
common, we repeat the argument to the left part of $x_{j-1}$ and
$y_{j-1}$.

(II) Thus, we get a symbolic system topologically conjugate to
$\mathcal Z_j$ with the same number of $j$-symbols such that the
modified sequences $x_{j-1}$ and $y_{j-1}$ have all $j$-cuts
common. On the other hand,  the sequences $x_{j-1}$ and $y_{j-1}$
remain $(j-1)$-compatible as we did not change the upper $j$-rows.
Furthermore, $x_{j-1}$ and $y_{j-1}$ remain $j$-separated. The
latter  means that there are two different $j$-symbols which are
$(j-1)$-compatible.

Now we produce a topological factor $\mathcal Z_j'$ of $\mathcal
Z_j$ by identifying  $j$-symbols which are $(j-1)$-compatible.
Thereby, we strictly reduce the number of $j$-symbols and
$j$-sequences from $\mathcal Z_j'$ have at most  $K-1$ different
$j$-symbols.

(III) It follows from (II) that the sequences  $x_{j-1}$ and
$y_{j-1}$ are not $j$-separated any more. The every  pair
$(x_i,y_i)$, $i\in [n_0,j-2]$ (in fact, their images in the factor
$\mathcal Z_j'$) remains $(i+1)$-separated since these rows have
not been changed. Moreover, the pair  $(x_i,y_i)$ still has a
common $j$-cut because in the construction we only added more
cuts.

We can now apply the same arguments as in (I) and (II) to the pair
$(x_{j-2},y_{j-2})$. As a result, we get a new factor with at most
$K-2$ $j$-symbols.

Repeating the same argument no more than $K-1$ times, we obtain a
new factor in which the pair $(x_{n_0},y_{n_0})$ remains
$(n_0+1)$-separated with a common $j$-cut while the family of
$j$-symbols consists of one element only. This is a contradiction.
\ep

\medskip
The next proposition is  called the ``Infection lemma'' in
\cite{downarowich_maass:2006}. Before proving the result, we
recall the notion of an eventually periodic sequence.

\mydefinition A sequence $x=\{x(n)\}_{n\in\mathbb Z}$ is called
{\it eventually periodic} if there are $n_0$ and $m$ such that
$x(n+m)=x(n)$ for all $n\geq n_0$. \medbreak

Recall also that the map $\pi_i:\mathcal Z_j\rightarrow \mathcal
Z_i$ denotes the projection to the first $i+1$ lines.

\begin{proposition}\cite{downarowich_maass:2006}\label{PropositionInfectionLemma}
Let $\{\mathcal A_0,\ldots,\mathcal A_j\}$ be  agreeable families of
$j$-symbols with $|\mathcal A_i|\leq K$, $i=0,\ldots,j,$ where
$j\geq K^2$. Suppose that the set $\{z_0,\ldots,z_{K^2}\}$
consists of  $i$-compatible ($i< j$) and pairwise $j$-separated
$j$-sequences with no common $j$-cuts. Then $\pi_i(z_0)$ is
eventually periodic.
\end{proposition}
\proc{Proof.} (a) We can set $\hat{z}:=\pi_i(z_k)$ because all the
points $z_k$ are $i$-compatible. Analogously, denote by $\hat v$
the restriction of a $j$-symbol $v$ to its top $i+1$ rows, i.e.
the projection of $v$ to $\mathcal A_i$. Draw a diagram $\mathcal
D$ consisting of the $j$-th rows of the elements $z_k$ one above
another with aligned zero coordinate.

{\small
$$\begin{array}{cccccccccccccccccccccccccccccc}
\hline \ldots & \multicolumn{6}{|c|} {z_0(-1)} &
\multicolumn{8}{|c|} {z_0(0)} & \multicolumn{7}{|c|} {z_0(1)} &
\multicolumn{5}{|c|} {z_0(2)} &
\ldots\\
\hline & \multicolumn{1}{c|}\ldots & \multicolumn{6}{|c|}
{z_1(-1)} & \multicolumn{5}{|c|} {z_1(0)} &
\multicolumn{6}{|c|}{z_1(1)} & \multicolumn{7}{|c|} {z_1(2)} &
\ldots
\\
 \hline & & & & & & & & & &  & & & \vdots & & \vdots & & & & & & & & & & & &
 \\ \hline  & & \multicolumn{1}{c|}\ldots &   \multicolumn{6}{|c|}
{ z_{K^2}(-1)} & \multicolumn{7}{|c|} {z_{K^2}(0)} &
\multicolumn{5}{|c|} {z_{K^2}(1)} & \multicolumn{4}{|c|}
{z_{K^2}(2)} & \ldots\\ \hline \\ \multicolumn{28}{c}{\mbox{
Diagram }\mathcal D\mbox{: } j\mbox{-th rows without common }
j\mbox{-cuts.}}
\end{array}$$ }

(b) Since we have $K^2+1$ sequences of $j$-symbols
$\{z_0,\ldots,z_{K^2}\}$, zero coordinates of
$\{z_0,\ldots,z_{K^2}\}$ are covered by at least $K+1$ copies of a
$j$-symbol $v$. Notice that since the elements
$\{z_0,\ldots,z_{K^2}\}$ have no common $j$-cuts, any of  two
copies of $v$ are shifted by some positive integer $0<l<|v|$.
Using the fact that $\pi_i(z_k)=\pi_i(z_0)$, we get that the
projection of the $j$-symbol $v$ to $\mathcal A_i$ satisfies the
``$l$-periodicity law'': $\hat v(n)=\hat v(n+l)$ for every $n\in
[0,|v|-1-l]$.

(c) Let $l_v$ be the minimal shift for $v$ appearing in the
diagram $\mathcal D$. This means that $\hat v(n)=\hat v(n+l_v)$
for every $n\in [0,|v|-1-l_v]$. As zero coordinate of $\hat z$ is
covered by at least two copies of $v$, we have that $\hat
z(0)=\hat z (l_v)$. Denote by $I$ the largest interval of $\mathbb
Z$ such that $0\in I$ and if $n\in I$, then $\hat z(n)=\hat
z(n+l_v)$.

(d) If $I$ is not bounded to the right, then the sequence $\hat
z=\pi_i(z_0)$ is eventually periodic and we are done.

If $I$ has the right end, set $m= (\max_{n\in I}n) +1$. Hence
$\hat z(m)\neq\hat z(m+l_v)$. Restrict the diagram $\mathcal D$ to
those $K+1$ elements $z_k$ in which the coordinate 0 is covered by
the $j$-symbol $v$. Since this diagram consists of $K+1$ lines,
the coordinate $m$ is covered by at least two copies of a
$j$-symbol $w$.

If $w=v$,  then by the choice of $l_v$, we get that $\hat
z(m)=\hat z (m+l_v)$, which is impossible.

Suppose that $w\neq v$. Let $r$ be the relative coordinate of $m$
within the extreme left copy of $w$ covering $m$ (see the
Figure\footnote{The idea of the figure is taken from
\cite{downarowich_maass:2006}.} below). Since there is no common
$j$-cuts, $r>0$. Hence, absolute coordinates of the extreme left
copy of $\hat w[0,r-1]$ within $\hat z$ intersects $I$. As $w\neq v$,
we have that the absolute coordinates of $\hat w[0,r-1]$ lie in $I$,
which implies that $\hat w(n)=\hat w(n+l_v)$ for all $n\in
[0,r-1]$.

{\small

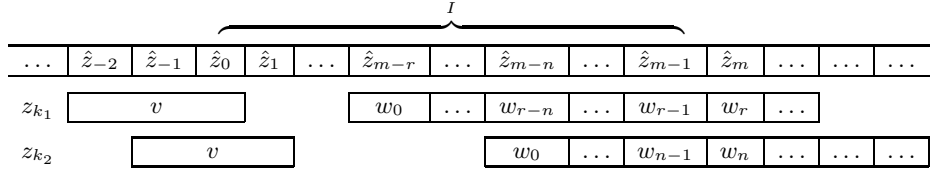
$$  \begin{array}{ccccccccccccccc} & & & \multicolumn{8}{c}{\overbrace{\qquad\qquad\qquad\qquad\qquad\qquad\qquad\qquad\qquad\quad}^I } \\
 \hline  \ldots & \multicolumn{1}{|c|}{\hat z_{-2}} &
\multicolumn{1}{|c|}{\hat z_{-1}} & \multicolumn{1}{|c|}{\hat z_0}
& \multicolumn{1}{|c|}{\hat z_1} & \ldots &
\multicolumn{1}{|c|}{\hat z_{m-r}} & \ldots &
\multicolumn{1}{|c|}{\hat z_{m-n}} & \ldots &
\multicolumn{1}{|c|}{\hat z_{m-1}} & \multicolumn{1}{|c|}{\hat
z_m} & \ldots & \multicolumn{1}{|c|}{\ldots}& \ldots
\\\hline
\\[-5pt]\cline{2-4}\cline{7-13} z_{k_1} & \multicolumn{3}{|c|} {v}
 & & & \multicolumn{1}{|c|} {w_0} & \ldots &  \multicolumn{1}{|c|} {w_{r-n}} & \multicolumn{1}{|c|} \ldots & \multicolumn{1}{|c|}
 {w_{r-1}} & \multicolumn{1}{|c|} {w_r} & \multicolumn{1}{|c|}\ldots
\\ \cline{2-4}\cline{7-13}
\\[-5pt] \cline{3-5}\cline{9-15} z_{k_2} & &
\multicolumn{3}{|c|}{v} &  & &  & \multicolumn{1}{|c|} {w_{0}} &
\multicolumn{1}{|c|} {\ldots} & \multicolumn{1}{|c|} {w_{n-1}} &
\multicolumn{1}{|c|} {w_n} & \multicolumn{1}{|c|} \ldots &
\multicolumn{1}{|c|} \ldots & \multicolumn{1}{|c|} \ldots
\\ \cline{3-5}\cline{9-15}

\end{array}$$ }

Now, consider the other copy of $w$. Since it is the extreme right
copy of $w$, the position $m$ is aligned with some relative
position $n\in [0,r-1]$, where the $l_v$-periodicity law holds.
Thus, $\hat z(m)=\hat z(m+l_v)$, which is a contradiction. \ep

%
%

\section{Finite Rank Systems}\label{SectionExpansiveFiniteRank}

In the section we generalize the result of
\cite{downarowich_maass:2006} to  homeomorphisms of a Cantor set
with  finite rank. To prove our main result of this section we develop the  method used in \cite{downarowich_maass:2006} to the case of aperiodic homeomorphisms. \medbreak

Let $B = (V,E, \leq)$ be an ordered Bratteli diagram such that the
ordering $\leq$ admits the continuous Vershik map $\varphi_B :
X_B\rightarrow X_B$ where $X_B$ is the space of infinite paths.
Denote  by $X_{\max}$ and $X_{\min}$ the sets of all maximal and
minimal paths of  $X_B$, respectively. Recall that by definition
of the Vershik map we have that $\varphi_B(X_{\max})=X_{\min}$.

Let $\{\mathcal A_j\}_{j\geq 0}$ be the family of agreeable
$j$-symbols associated to the diagram $B$ (see the details in Example
\ref{j-symbols in Bratteli diagrams}).
\medskip

\myremark (1) Any infinite path $x=(x_n)\in X_B\setminus Orb_{\f_B}(X_{min}\cup
X_{max})$ can be represented as an infinite matrix denoted by
$[x]$ with rows indexed by $i$ from $0$ to $\infty$ and columns
indexed by $j\in (-\infty,\infty)$. The matrix $[x]$ is formed by
an increasing sequence of $i$-symbols ($i=0,1,\ldots$) which
cross the coordinate $0$.  The $i$-symbol  corresponds to the
vertex $v\in V_i$ which is crossed by the path $x$ at the level
$i$, and the position of $0$ coordinate in this $i$-symbol is
defined by the order of the finite path $(x_1,\ldots,x_i)$ amongst
all paths connecting $v_0$ and $v$. We see that this construction
gives not only a single point $x\in X_B$ but the entire orbit
$Orb_{\f_B}(x) = \{\f_B^i(x)\}_i$.

(2) We notice  that for every $x\in X_{\max}$ the Vershik map
uniquely defines $y \in X_{\min}$ such that $y= \f_B(x)$. If we
applied the construction used in (1) to the paths $x$ and $y$, we
would get two one-sided matrices   $[y]^+$ and $[x]^-$ infinite to
the right and left, respectively. As $Orb_{\f_B}(x)
=Orb_{\f_B}(y)$ it is natural to assign the concatenated matrix
$[x]^-[y]^+$ to $x$ (or $y$).

(3) Let $X$ be the set of infinite matrices described in (1) and
(2). Then the Bratteli-Vershik model $(X_D,\f_B)$ is conjugate to
$(X, T)$ where $T$ is the left shift in $X$. This fact allows us to identify the sets $X$ and $X_B$. \medbreak

Recall that $\mathcal Z_i$ denotes the set of all $i$-sequences,
see Definition \ref{Definition_J-sequences}. Let $\pi_i:
X_B\rightarrow \mathcal Z_i $ be the map which restricts each
matrix $[x]$, $x\in X_B$, to the first $i+1$ rows. Set
$X_i=\pi_i(X_B)$. Clearly, $X_i$ is a closed shift-invariant
subset of $\mathcal Z_i$.   Observe also that $(X_i,T)$ is a
factor of $(X_B,\varphi_B)$.\medbreak

\mydefinition Let  $d$ be a metric on $X$ which generates the
topology. It is said that a homeomorphism $S: X\rightarrow X$ is
{\it expansive} if there exists $\delta>0$ such that for any
distinct $x,y\in X$ there is $m\in\mathbb Z$ with
$d(S^mx,S^my)>\delta$. The number $\delta$ is called an {\it
expansive constant}. \medbreak

Note that the notion of expansiveness does not depend on the
choice of the metric $d$, see \cite[Section 5.6]{walters:book}
\medbreak

\myremark   If $(X_B,\f_B)$ is an expansive system, then, due to
the famous theorem of Hedlund,  $(X_B,\f_B)$ is homeomorphic to
$(X_i,T)$ for all sufficiently large $i$, for the details see the
proof of Theorem 5.24 in \cite{walters:book}. \medbreak

The following definitions agree with
Definition \ref{Definition_Compatible_J-sequences}.

\mydefinition  We say that two distinct points $x$ and $y$ from $X_B$ are {\it $i$-compatible} if
$\pi_i(x)=\pi_i(y)$. If $\pi_i(x)\neq \pi_i(y)$, then $x$ and $y$
are called {\it i-separated}. Clearly, any distinct points $x$ and
$y$ are $i$-separated for some $i$. The largest integer $i$ such
that $\pi_i(x)=\pi_i(y)$ is called the {\it depth} of $x$ and $y$.
A pair $(x,y)$ has a common {\it $j$-cut} if there exist
$n\in \mathbb Z$ and $j$-symbols $v$ and $w$ such that $v$ appears
at the position $n$ in $\pi_j(x)$ and $w$ appears at $n$ in
$\pi_j(y)$.  Notice that if $x$ and $y$ have a common $j$-cut,
then $x$ and $y$ have a common $j'$-cut for all $j'\leq j$.
\medbreak

\mydefinition A Cantor dynamical system $(X,S)$ has the {\it
topological rank} $K>0$ if it admits a Bratteli-Vershik model
$(X_B,\varphi_B)$ such that the number of vertices of the diagram
$B$ at each level is not greater than $K$ and $K$ is the least
possible number of vertices for any Bratteli-Vershik realization.

Clearly, if a system $(X,S)$ has the rank $K$, then, by an
appropriate telescoping, we can assume that the diagram $B$ has
exactly $K$ vertices at each level. \medbreak

The next statement shows that the number of minimal
components of a finite rank system is bounded (see also Proposition
\ref{PropositionNumberMinimalComponentsGeneralSubst} for a similar
result for substitutional systems).

\begin{proposition}\label{PropositionNumberOfMinimalComponentsVershikMap} Let $(X,S)$ be a Cantor aperiodic dynamical system of a
finite rank $K$. Then $(X,S)$ has at most $K$ minimal components.
\end{proposition}
\proc{Proof.} By definition,  $(X,S)$ can be realized as a
Bratteli-Vershik model with at most $K$ vertices at each level.
 Assume that there exist $(K+1)$-minimal components $Z_0,\ldots, Z_{K}$ for a homeomorphism $S$.
Then for each level $i$ there exist a vertex $v_i\in V_i$ and two
paths $x_i$ and $y_i$ from $X_B$ such that they belong to
different minimal components from $Z_0,\ldots, Z_K$ and pass
through $v_i$. It follows that there exist distinct $n$ and $m$
such that $0\leq n, m \leq K$ and $x_i \in Z_n, \ y_i\in Z_m$ for
infinitely many indexes $i\in I$.

Define the metric $d$ on the space $X_B$ as follows:
\begin{equation}\label{Equation_Metric}
d(\{x_n\},\{y_n\})=\frac{1}{\min\{k : x_k\neq y_k \}}.
\end{equation}

Let $i_0\in I$ be chosen such that $\mbox{dist}(Z_n,Z_m) > 1/i_0$.
Notice that we can find a path $z\in Orb_S(x_{i_0})$ such that the
first $i_0$ edges of $z$ coincide with those of $y_{i_0}$. It
follows that $\mbox{dist}(z, Z_m) \leq 1/i_0$, which is a
contradiction. \ep
\medskip

We recall the definition of the enveloping semigroup, see, for
example, the book \cite{glasner:joinings} for a wider coverage of
the subject.

\mydefinition \label{DefinitionEnvelopingSemigroup} Let $S: X\to X$ be a homeomorphism. By definition, the {\it enveloping semigroup}
$E=E(X,S)$ of the dynamical system $(X,S)$ is the closure of the set $\{S^n: n\in\mathbb Z\}$ in $X^X$ with respect to the topology of
pointwise convergence.
\medbreak

The main result of the section is the following statement.

\begin{theorem}\label{TheoremFiniteRank} Let $(X,S)$  be an aperiodic Cantor dynamical system of finite rank $K$. If the restriction of $(X,S)$ to every minimal component is not homeomorphic to an odometer, then $(X,S)$ is expansive.
\end{theorem}
\proc{Proof.} Without loss of generality, we can assume  that $X$
is the path-space of an ordered Bratteli diagram $B=(V,E,\leq)$
and $S$ is the Vershik map defined by the ordering $\leq$.

Assume that the system $(X,S)$ is not expansive. Then for each
$i\geq 1$ there exists a pair $(x_i,y_i)$ such that
$d(S^mx_i,S^my_i)\leq 1/i$ for all $m\in\mathbb Z$, where the
metric is defined  in Equation \ref{Equation_Metric}. This is
equivalent to the fact that the pair $(x_i,y_i)$ is
$i$-compatible. Therefore, for infinitely many $i$ there is a pair
of depth $i$. By telescoping of $B$, can  assume that for each
$i\geq 1$ there is a pair  of depth $i$.

Now we have two opposite statements:

{\it (1) There exists $i_0$ such that for all $i\geq i_0$ and
every $j\geq i$ there is a pair of depth $i$ with a common
$j$-cut.

(2) For all $i_0$ there exist $i_0\leq i<j$ such that any pair of
depth $i$ has no common $j$-cuts.}

\medskip  It follows from Proposition
\ref{PropositionSymbolicFormasimNoCommonCuts} that statement (1) is
never true. Thus, for all $i_0$ there exist $j> i \geq i_0$ such
that any pair of depth $i$ has no common $j$-cuts. Therefore, any
pair of depth $i$ has no common $j'$-cuts for $j'>j$.

After telescoping of the diagram, we can assume that any pair of
depth $i$ has no common $(i+1)$-cuts. Recall that by our
assumption of non-expansiveness of $S$, there exists at least one
pair of depth $i$.

Let $E=E(X,S)$ be the enveloping semigroup for $(X,S)$. Therefore,
$Ex=\{\gamma(x) : \gamma\in E\}=\overline{Orb_S(x)}$. It follows
that $Ex$ contains at least one minimal $S$-component. By Proposition
\ref{PropositionNumberOfMinimalComponentsVershikMap} there exist at most $K$ minimal $S$-components,  $Z_0,\ldots,
Z_{K-1}$. Fix any path $z_i\in Z_i$, $i=0,\ldots, K-1$.

Fix any level $i_0$. For every $i\in [i_0,i_0+K^2+K]$ take a pair
$(x_i,y_i)$ of depth $i$ without common $(i+1)$-cuts. For each
pair $(x_i,y_i)$ find an element $\gamma_i\in E$ such that
$\gamma_i(x_i)=z_j$ for some $z_j$, $j=j(i)$.

(a) Observe that the pair $(\gamma_i(x_i),\gamma_i(y_i))$ has the
same depth as $(x_i,y_i)$ and has no common $(i+1)$-cuts. Indeed,
$\gamma_i$ is the pointwise limit of a sequence $S^{n_k}$. As each
projection $\pi_m: X\rightarrow X_m$ is continuous, we get
$$\begin{array}{ll}\pi_i(\gamma_i(x_i))=\pi_i(\lim S^{n_k}x_i)=\lim
\pi_i(S^{n_k}x_i)=\lim T^{n_k}\pi_i(x_i)=\\
=\lim T^{n_k}\pi_i(y_i)=\lim\pi_i(S^{n_k}y_i)=
\pi_i(\gamma_i(y_i)).\end{array}$$ This, in particular, implies
that $(\gamma_i(x_i),\gamma_i(y_i))$ is $i$-compatible.

It remains to show only that $(\gamma_i(x_i),\gamma_i(y_i))$ has
no common $(i+1)$-cuts. Indeed, suppose $\gamma_i(x_i)$ and
$\gamma_i(y_i)$ have a common $(i+1)$-cut. Then
$T^{n_k}\pi_{i+1}(x_i)$ and $T^{n_k}\pi_{i+1}(y_i)$ have a common
$(i+1)$-cut for all $k$ sufficiently large. This implies that the pair
$(x_i,y_i)$ has a common $(i+1)$-cut, a contradiction.

(b) Consider the pairs $(\gamma_i(x_i),\gamma_i(y_i))$, $i\in
[i_0,i_0+K^2+K]$. As each $\gamma_i(x_i)$ is equal to some $z_j$,
$j=j(i)$,  and we have at most $K$ of them, we can  choose $K^2$
pairs $(\gamma_i(x_i),\gamma_i(y_i))$, $i\in I$, such that
$\gamma_i(x_i)$'s are the same for all $i\in I$, $|I|=K^2$, say
$\gamma_i(x_i)=z_{k_0}\in Z_{k_0}$. Observe that the pair
$(\gamma_i(y_i),\gamma_j(y_j))$, $i,j\in I$, $i<j$ has the depth
$i$ and has no common $j$-cuts.

(c) Setting $z_i'=\gamma_i(y_i)$ for $i\in I$, we obtain a family
$C=\{z_i' : i\in I\}\cup\{z_{k_0}\}$ of $i_0$-compatible and
pairwise $(i_0+K^2 + K)$-separated elements with no common
$(i_0+K^2+K)$-cuts. Considering the projection $\pi_{i_0+K^2 + K}$
of $C$ and applying Proposition \ref{PropositionInfectionLemma},
we get that $\pi_{i_0}(z_{k_0})$ is eventually periodic, and, by
minimality, $\pi_{i_0}(Z_0)$ is periodic.

\medskip Thus, we get that for each $i$, there is a minimal
$S$-component $Z_i$ such that $\pi_i(Z_i)$ is $T$-periodic. As we
have at most $K$ minimal $S$-components, there is a minimal
component, say $Z_{k_0}$, such that $\pi_i(Z_{k_0})$ is periodic
for infinitely many $i$. This implies that the restriction of $S$
to $Z_{k_0}$ is homeomorphic to an odometer. The theorem is proved.
\ep

%
%

\section{Recognizability of aperiodic substitutions}\label{SectionNotionOfSubstDynamSystem}
In the section we study  dynamical properties of an arbitrary
substitutional dynamical system $(X_\sigma,T_\sigma)$. First of
all, we estimate the number of minimal components of the system
$(X_\sigma,T_\sigma)$. Then we build a sequence of K-R partitions
of $X_\sigma$ using geometrical properties of the substitutional dynamical system.

\medskip Let $A$ denote a finite alphabet and $A^+$  the set
of all non-empty words over $A$. Set also $A^*=A^+\cup \{\emptyset\}$. For a word $w=w_0\ldots w_{n-1}$
with $w_i\in A$ let $|w|=n$ stand for its length. For any two
words $v,w\in A^+$, the symbol `$v\prec w$' means that $v$ is a
factor of $w$.

\mydefinition By a {\it substitution} we mean any map $\sigma :
A\rightarrow A^+$.\medbreak

Any  map $\sigma : A\rightarrow A^+$ is extended to the map
$\sigma :A^+\rightarrow A^+$ by concatenation. We define the {\it
language} $L(\sigma)$   of  a substitution $\sigma$ as the set of
all words which appear as factors of $\sigma^n(a)$, $a\in A$, $n\geq 1$. By definition, we also set that $\sigma^0(a)=a$ for all
$a\in A$.

For any substitution $\sigma : A\rightarrow A^+$, define
\begin{equation}\label{EquationMaxMinSubstitution} |\sigma|=\min\limits_{a\in
A}|\sigma(a)|,\qquad ||\sigma||=\max\limits_{a\in
A}|\sigma(a)|.
\end{equation}
Observe that the functions $n\mapsto |\sigma^n|$ and $n\mapsto
||\sigma^n||$ are not decreasing.

\mydefinition By a {\it substitutional dynamical system associated
to a substitution $\sigma$}, we mean a pair $(X_\sigma,T_\sigma)$,
where
$$X_\sigma=\{x\in A^\mathbb Z : x{[-n,n]}\in L(\sigma)\mbox{ for
any }n\}$$ and $T_\sigma$ is the shift on $A^\mathbb {Z}$. We will
denote the $k$-th coordinate of  $x\in X_\sigma$ by $x[k]$ or by
$x(k)$. For any words $v,w\in L(\sigma)$, set $[v.w]:= \{x\in
X_\sigma : x[-|v|,|w|-1]=vw\}$. \medbreak

\mydefinition \label{DefinitionAperiodicSubstitutions} A
substitution $\sigma : A\rightarrow A^+$ is called {\it aperiodic}
if the system $(X_\sigma,T_\sigma)$ has no periodic points.
\medbreak

\myremark (1) In general, the set $X_\sigma$ can be empty. To
avoid trivialities, {\it we will always assume that the
substitution $\sigma$ is such that  $X_\sigma$ is an uncountable
subset of $A^{\Z}$}. In particular, we always have that
$||\sigma^n||$ tends to the infinity.

(2) If $L(\sigma^k)=L(\sigma)$ for some $k\geq 1$, then the
substitutional dynamical systems associated to $\sigma$ and
$\sigma^k$ coincide, i.e., $X_\sigma=X_{\sigma^k}$.

(3) The set $X_\sigma$ is
a 0-dimensional compact metrizable space whose topology is generated by the metric
\begin{equation}d(\{x[n]\},\{y[n]\})= \sum_{n=-\infty}^\infty
\frac{d'(x[n],y[n])}{2^{|n|}}
\end{equation}
where $d'$ is the discrete metric on $A$, i.e. $d'(a,b) = 1$ iff $a\neq b$. Notice that the set
$X_\sigma$ can have isolated points. We do not require
$X_\sigma$ to be perfect.\medbreak

Let $\sigma :A\rightarrow A^+$ be a substitution. We denote by
$A_l$ the set of all letters $a\in A$ such that
$|\sigma^n(a)|\to\infty$ as $n\to\infty$. Set also $A_s=A\setminus
A_l$. Here the subindexes `s' and `l' stand for `short' and
`long', respectively. Observe that $\sigma(a)\in A^+_s$ for all
$a\in A_s$.

To estimate the number of minimal components of an arbitrary substitutional system, we need the
following proposition which asserts that the length of words formed by short letters is  uniformly  bounded. This result will be also applied to Bratteli diagram construction (Theorem \ref{TheoremFromSubtitutionsToDiagrams}).

\begin{proposition}\label{Proposition_BoundednceOfShortWords} Let $\sigma: A\rightarrow A^+$ be an aperiodic
substitution. Then there exists $M>0$  such that every word $W\in
L(X_\sigma)$ with $|W|\geq M$ contains at least one letter from
$A_l$.
\end{proposition}
\proc{Proof.} Assume the converse, i.e., for any $m>0$  there is a
word $W_m\in L(X_\sigma)\cap A_s^+$ with $|W_m|=m$.

(0) Since $|\sigma^n(a)|\leq K<\infty$ for all $a\in A_s$ and some $K$, there are
positive integers $n_a$ and $p_a$ such that
$\sigma^{n_a}(a)=\sigma^{n_a+kp_a}(a)$ for all $k\geq 0$, $a\in
A_s$. Set $n_0=\max\{n_a : a\in A_s\}$. Then it follows that
$$\sigma^{n}(a)=\sigma^{n+kp}(a)\mbox{ for all }k\geq 0\mbox{ and }a\in A_s$$ where
$n\geq n_0$ and $p=\prod_{a\in A_s} p_a$. Setting
$\tau=\sigma^{pn_0}$, we obtain that
$\tau^2(a)=\sigma^{pn_0+pn_0}(a)=\tau(a)$ for all $a\in A_s$, i.e.
$\tau^k(a)=\tau(a)$ for all $a\in A_s$ and $k\geq 1$.

(1) Since $|A| < \infty$, we find a letter $a\in A_l$ and an
infinite set $I\subseteq \mathbb N$ such that
$W_m\prec\sigma^{k_m}(a)$ for $m\in I$. Write down each word
$\sigma^n(a)$ as
$$\sigma^n(a)=L_1^{'(n)}S_1^{'(n)}\ldots L_k^{'(n)}S_k^{'(n)}L_{k+1}^{'(n)}\ldots
L_{r_n}^{'(n)}S_{r_n}^{'(n)}$$ where $L^{'(n)}_i\in A_l^+$ and
$S_i^{'(n)}\in A_s^+$. Observe also that the words $L_1^{'(n)}$
and $S_{r_n}^{'(n)}$ can be empty. Since $\sigma(S_i^{'(n)})\in
A_s^+$, the maximal length of
$\{S_1^{'(n)},\ldots,S_{r_n}^{'(n)}\}$ must monotonically  tend to
infinity as $n\to\infty$.

Let $\tau=\sigma^{n_0p}$, where the integers $n_0$  and $p$ are
defined above. Then, we decompose each word $\tau(a)$, $a\in A_l$,
into blocks over short and long letters
\begin{equation}\label{tau-decomposition}
\tau^n(a)=L_1^{(n)}S_1^{(n)}\ldots
L_k^{(n)}S_k^{(n)}L_{k+1}^{(n)}\ldots L_{d_n}^{(n)}S_{d_n}^{(n)}
\end{equation}
where $L^{(n)}_i\in A_l^+$ and $S_i^{(n)}\in A_s^+$. We notice
that the maximal length of $\{S_1^{(n)},\ldots,S_{d_n}^{(n)}\}$
tends to infinity as $n\to\infty$.

Since $\tau(S_k^{(n)})\in A_s^+$, each block $S_k^{(n)}$ is {\it
followed by a block} $S_i^{(n+1)}$ (in symbols,
$S_k^{(n)}<S_i^{(n+1)}$) in the following sense: let $a_k^{(n)}$
and $b_k^{(n)}$ denote the last and the first letters of
$L_k^{(n)}$, respectively, then $S_i^{(n+1)}$ is the maximal
factor of
$\tau(a_k^{(n)}S_k^{(n)}b_{k+1}^{(n)})=\tau(a_k^{(n)})\tau(S_k^{(n)})
\tau(b_{k+1}^{(n)})$ over $A_s^+$ that contains 
$\tau(S_k^{(n)})$. In other words, $S_i^{(n+1)}$ is the maximal
factor of $\tau^{n+1}(a)$ over $A^+_s$ that appears between
$a_i^{(n+1)}$ and $b_{i+1}^{(n+1)}$.

(2) For $N$ sufficiently large, there exists $1\leq k_0\leq d_{N}$
(see (\ref{tau-decomposition}) for the definition of $d_N$) such
that the block $S_{k_0}^{(N)}$ is followed by blocks
$S_{k_i}^{(N+i)}$,
$$S_{k_0}^{(N)}<S_{k_1}^{(N+1)}<S_{k_2}^{(N+2)}<\ldots$$ with
$|S_{k_i}^{(N+i)}|\to\infty$ as $i\to\infty$.

Therefore,  we see from the above argument that
$$\begin{array}{cc}\tau^i(a_{k_0}^{(N)})=\boxed{*} a_{k_i}^{(N+i)}V_i,\\
\tau^i(b_{k_0+1}^{(N)})=Z_i
b_{k_i+1}^{(N+i)}\boxed{**}\end{array}$$ where $V_i,Z_i\in A_s^*$
and $\boxed{*}$, $\boxed{**}$ are some words. Then we can  find
positive integers $q$ and $t$ such that
$$\begin{array}{cc}a_{k_t}^{(N+t)}=a_{k_{t+ql}}^{(N+t+ql)},\\
b_{k_t+1}^{(N+t)}=b_{k_{t+ql}+1}^{(N+t+ql)}\end{array}$$  for all
$l\geq 0$.

(3) Since $|S_{k_i}^{(N+i)}|\to\infty$, we have that $a_{k_i}^{(N+i)}$ or $b_{k_i+1}^{(N+i)}$ adds to
$S_{k_0}^{(N)}$ new blocks over $A_s^+$, i.e., $|V_i|\to\infty$ or $|Z_i|\to\infty$ as
$i\to\infty$. For  definiteness, we  assume that so does
$a_{k_i}^{(N+i)}$. Therefore,
$$\tau^{q}(a_{k_t}^{(N+t)})=
\boxed{*}a_{k_{t+q}}^{(N+t+q)}Q=\boxed{*}a_{k_t}^{(N+t)}Q$$
and, inductively,
$$\tau^{lq}(a_{k_t}^{(N+t)})=
\boxed{**}a_{k_{t}}^{(N+t)}Q\tau^q(Q)\ldots\tau^{q(l-1)}(Q),\;l\geq
1,
$$
where $Q\in A_s^+$ and $\boxed{*}$, $\boxed{**}$ are some words.

It follows from the definition of $\tau$ that $\tau^l(Q)=\tau(Q)$
for every $l\geq 1$. Therefore, $X_\sigma$ contains a periodic
sequence $\tau(Q)^\infty$, which is impossible. \ep

\begin{proposition}\label{PropositionNumberMinimalComponentsGeneralSubst}
 Let $\sigma : A\rightarrow A^+$ be an aperiodic substitution.
Then $(X_\sigma,T_\sigma)$ has no more than $|A|$ minimal
components.
\end{proposition}
\proc{Proof.} Assume the converse. Set $K=|A|$. Take any $K+1$
minimal components $Z_0,\ldots,Z_K$ and let $z_i\in Z_i$,
$i=0,\ldots,K$.

 We claim that there exist an infinite set $I\subseteq\mathbb N$ and 
letters $b_j\in A_l$, $j=0,\ldots,K$, such that $\sigma^i(b_j)\prec
z_j$  for all $i\in I$, $j=0,\ldots,K$.

Indeed, by Proposition \ref{Proposition_BoundednceOfShortWords},  find $M>0$ such that any word $w\in L(X_\sigma)$ of length at least $M$
contains a letter from $A_l$.   Given $n>0$, take $m>M||\sigma^n||$. By definition of $X_\sigma$, $z_0[-m,m]$ is a factor of $\sigma^k(a)$
for some $a\in A$ and $k$. Since the function $n\mapsto
||\sigma^n||$ is not decreasing, we get that $k>n$. Let
$\sigma^{k-n}(a)=a_0\ldots a_{d-1}$ with $a_i\in A$. We can choose the maximal interval $[i,j]$ of $[0,d-1]$ such that
$\sigma^n(a_i\ldots a_j)$ appears in $z_0[-m,m]$. The choice of $m$ guarantees us that $j-i\geq M$. Therefore, at least one of the letters $a_r$ belongs to $A_l$, $i\leq r\leq j$.

Thus, for all $n\in\mathbb N$
there is $a_n\in A_l$ such that $\sigma^n(a_n)$ appears in $z_0$. Therefore, there is an infinite set $I_0\subseteq \mathbb N$ and a letter $b_0\in A_l$ such that
$\sigma^i(b_0)$ appears in $z_0$ for all $i\in I_0$. Analogously, for all
$i\in I_0$ there is $a_i\in A_l$ such that $\sigma^i(a_i)$ appears in
$z_1$. Therefore, there is an infinite set $I_2\subseteq I_1$ and a letter $b_1\in A_l$ such that $\sigma^i(b_1)$ appears in $z_1$ for all $i\in I_1$. Repeating the argument, we find  $I=I_K$ and letters $b_0,\ldots,b_{K-1}$. This proves
the claim.

Since we have only $K$ letters, there is a letter, say $a\in A_l$,
such that $\sigma^i(a)$ appears in two distinct sequences $z_l$
and $z_d$ for all $i\in I$. Fix any $i\in I$. Denote by $E[i]$ the
integer part of $|\sigma^i(a)|/2$. Note that $|\sigma^i(a)|\to\infty$. By shifting $z_l$ and $z_d$, if
necessary, we can assume that $\sigma^i(a)$ appears in $z_l$ and
$z_d$ in such a way that the $E[i]$-th coordinate of $\sigma^i(a)$
is aligned with   zero coordinates of $z_l$ and $z_d$. Therefore,
the $d$-distance between the compact sets $Z_l$ and $Z_d$ is less
than $1/{2^{E[i]}}$, for any $i\in I$, which is impossible. \ep

\myremark It immediately follows from the proof of Proposition \ref{PropositionNumberOfMinimalComponentsVershikMap} that the result is still true if we  replace the  aperiodicity of $\sigma$ by the condition that  $|\sigma^n|\to\infty$ as $n\to\infty$.

\mydefinition Let $\sigma: A\rightarrow A^+$ be a substitution,
$a\in A$, $k>0$. The words of the form $\sigma^k(a)$ will be called {\it $k$-words}. Let $n>k$ and
$\sigma^{n-k}(a)=a_0\ldots a_{m-1}$ with $a_i\in A$. Then
$\sigma^n(a)=\sigma^k(a_0)\ldots\sigma^k(a_{m-1})$. We will say
that the $n$-word $\sigma^n(a)$ is {\it naturally decomposed into
$k$-words} $[\sigma^k(a_0),\ldots,\sigma^k(a_{m-1})]$. \medbreak

\myremark\label{RemarkNaturalDecomposition} Suppose that an
$n$-word $\sigma^n(a)$ is naturally decomposed into $k$-words
$[\sigma^k(a_0),\ldots,\sigma^k(a_{m-1})]$ and into $k+1$-words
$[\sigma^{k+1}(c_0),\ldots,\sigma^{k+1}(c_{l-1})]$, i.e.,
$$\sigma^n(a)=\sigma^k(\sigma(c_0))\ldots\sigma^k(\sigma(c_{l-1}))=
\sigma^k(a_0)\ldots\sigma^k(a_{m-1}).$$

Since the decompositions are natural, we get that each letter
$a_i$ appears as a factor of some $\sigma(c_j)$. \medbreak

\medskip The following  result shows that each element of $X_\sigma$ can be written as a concatenation of $1$-words.

\begin{proposition}\label{Proposition-1-cutting} Let $\sigma : A\rightarrow A^+$ be a substitution.
Then for every $x\in X_\sigma$ there exist $y\in X_\sigma$ and
$i\in\{0,\ldots,|\sigma(y[0])|-1\}$ such that
$x=T^i_\sigma\sigma(y)$.
\end{proposition}
\proc{Proof.} {\it I.} Take any $x\in X_\sigma$. For each $n\geq
1$, find $m_n>n||\sigma^n||$. By definition of $X_\sigma$, there
are $v^{(n)}\in A$ and $k_n>0$ such that $x{[-m_n,m_n]}$ is a
factor of $\sigma^{k_n}(v^{(n)})$. It is evident that $k_n>n$. Let
$\sigma^{k_n-n}(v^{(n)})=v_1^{(n)}\ldots v_{d_n}^{(n)}$ with
$v_i^{(n)}\in A$. Hence, the word $x{[-m_n,m_n]}$ is a factor of
the word $\sigma^n(v_1^{(n)})\ldots\sigma^n(v_{d_n}^{(n)})=
\sigma^{k_n}(v^{(n)})$ which is considered as the concatenation of
$n$-words. Take the maximal interval $[i_n,j_n]\subseteq
\{1,\ldots,d_n\}$ so that the word
$\sigma^n(v_{i_n}^{(n)})\sigma^{n}(v_{i_n+1}^{(n)})
\ldots\sigma^n(v_{j_n}^{(n)})$ appears as a factor of
$x{[-m_n,m_n]}$. That is
$$x{[-l_n,|\sigma^n(v_{i_n}^{(n)}\ldots v_{j_n}^{(n)})|-1-l_n]}=
 \sigma^n(v_{i_n}^{(n)}\ldots v_{j_n}^{(n)}),$$
where $l_n$ is defined as the natural position of occurrence of
$x[0]$ within the word
$\sigma^n(v_{i_n}^{(n)})\ldots\sigma^n(v_{j_n}^{(n)})$.  Setting
$w_n=v_{i_n}^{(n)}\ldots v_{j_n}^{(n)}$, we get that
\begin{equation}\label{Claim1_Equation}x{[-l_n,|\sigma^n(w_n)|-1-l_n]}=
\sigma^n(w_n)\mbox{ for all }n\geq 1
\end{equation}
and $l_n\to\infty$, $(|\sigma^n(w_n)|-l_n)\to\infty$ as
$n\to\infty$.

{\it II.} For each $n$ and $i=i_n,\ldots,j_n$, take the natural
decomposition of the $n$-word $\sigma^n(v_i^{(n)})$ into 1-words
and write down all of these 1-words from left to right, say
$$\sigma^n(w_n)=\sigma^n(v_{i_n}^{(n)})\ldots \sigma^n(v_{j_n}^{(n)})=
\sigma(y^{(n)}_{-s_n})\ldots \sigma(y^{(n)}_{d_n}).$$
Enumerate these 1-words  by the following rule: as all of them
naturally appear in the sequence $x$, we set $\sigma(y^{(n)}_0)$
to be the 1-word that naturally covers $x[0]$; $\sigma(y_{-s}^{(n)})$
is the $s$-th 1-word on the left of $\sigma(y_0^{(n)})$;
and $\sigma(y_s^{(n)})$ is the $s$-th 1-word on the right of
$\sigma(y_0^{(n)})$ if $s>0$.

For every $n\geq 1$, define the
integer $f^{(n)}_{0}\leq 0$  as the position  of the natural occurrence
of 1-word $\sigma(y_{0}^{(n)})$  within $x$, i.e., $$x[f^{(n)}_0,|\sigma(y_{0}^{(n)})|-1+f_0^{(n)}]=\sigma(y_0^{(n)}).$$

{\it III.} Considering 1-words $\sigma(y_0^{(n)})$  and the positions of their occurrences in $x$,
we find an infinite subset $I_0\subseteq \mathbb N$ such that
$$y_0^{(n)}=y_0^{(n')}\mbox{ and }
f_0^{(n)}=f_0^{(n')}\mbox{ for
any }n,n'\in I_0.$$

Then by induction on $k$, we find an infinite set $I_k\subseteq
I_{k-1}\subseteq\ldots\subseteq I_0$ such that
$$y_i^{(n)}=y_i^{(n')}\mbox{ for any }n,n'\in I_k\mbox{ and
}i=-k,\ldots,k.$$

This defines a two-sided sequence
$$y=\{\ldots,y_{-k}^{(n_k)},y_{-k+1}^{(n_{k-1})},\ldots,y_0^{(n_0)},
y_1^{(n_1)},\ldots,y_{k}^{(n_{k})},\ldots
\}\mbox{ with }n_k\in I_k. $$ Observe that $y\in X_\sigma$ and
$x=T^i_\sigma\sigma(y)$, where $i=-f_0^{(n_0)}$. \ep

$$\begin{array}{cc|c|c|c|c|c|c|c|c}\hline \ldots & \multicolumn{1}{|c|} {x_{-3}} &
x_{-2}& x_{-1} & x_0 & x_1 & x_2 & x_3 & x_4  & \ldots \\
\hline \ldots & \multicolumn{2}{|c|} {\sigma(y_{-1})} &
\multicolumn{3}{|c|} {\sigma(y_0)} & \multicolumn{3}{|c|}
{\sigma(y_1)} & \ldots\\ \hline \\ \multicolumn{10}{c}{
\mbox{Decomposition into 1-words}}\\
\end{array}$$

\myremark\label{Remark_1-cutting} As an immediate corollary  of Proposition  \ref{Proposition-1-cutting} we get that for any $n>0$ and any $x\in X_\sigma$ there exists $y\in X_\sigma$ and $i\in \{0,\ldots,|\sigma^n(y[0])|-1\}$ such that $x=T^i_\sigma\sigma^n(y)$.

%
%

\mydefinition We say that a substitution is {\it recognizable} if
for each $x\in X_\sigma$ there exist a unique $y\in X_\sigma$ and
unique $i\in \{0,\ldots,|\sigma(y[0])|-1\}$ such that
\begin{equation}\label{Equation_Cutting_into_1-words}x=
T^i_\sigma\sigma(y).\end{equation}
\medbreak

\myremark We note that in the theory of primitive substitutions it
was a long-standing problem to establish the uniqueness of this
representation. It follows from the works \cite{host:1986} or
\cite{queffelec:book} that under the assumption of {\it bilateral
recognizability} (see \cite{mosse:1992}) and injectivity  of the
substitution on the alphabet, each aperiodic primitive
substitution is recognizable in our sense. However, due to the
works of Moss\'e \cite{mosse:1992} and \cite{mosse:1996} it became
clear that, in fact, each aperiodic primitive substitution is
recognizable. \medbreak

 Now we will show that an arbitrary aperiodic substitution
is recognizable. Our proof involves  the usage  of
Downarowicz-Maass'  techniques developed in Section
\ref{SectionSymbolicFormalism}. So, first of all, we introduce
agreeable families of $j$-symbols.

\mydefinition Let $\sigma : A\rightarrow A^+$ be a substitution.
Set $A_{i-1}=A$ and $\sigma_i=\sigma$, for all $i\geq 1$. Denote
by $\mathcal A_j$ the family of $j$-symbols determined by the
alphabets $\{A_i\}_{i\geq 0}$ and the maps $\{\sigma_i\}_{i\geq
1}$ as in Definition \ref{Definition_Agreeable_J-Symbols}.\medbreak

To make the explanation more comprehensible, we will label the boxes
from a row $i$ of $j$-symbol $[a]_j\in\mathcal A_j$ by the symbols
$\sigma^i(b)$, $b\in A$, instead of just $b\in A$. \medbreak

 \begin{center}{\small $$\begin{array}{cc|c|c|c|c|c|c|c|}\hline \multicolumn{1}{|c|} {b_0} & b_1 & b_2 & b_3 & \ldots & b_{m-3} & b_{m-2} & b_{m-1}\\
 \hline \multicolumn{2}{|c|} {\sigma(c_0)} & \multicolumn{2}{|c|} {\sigma(c_1)} &
 \ldots&
 \multicolumn{3}{|c|} {\sigma(c_{d-1})} \\
 \hline \multicolumn{8}{|c|} \vdots
 \\
 \hline \multicolumn{8}{|c|} {\sigma^j(a)}\\
 \hline \\   \multicolumn{8}{c} {j\mbox{-symbol}}\end{array}$$}
 \end{center}

 \medskip Denote by  $\mathcal Z_j$ the set of all $j$-sequences, see Definition
\ref{Definition_J-sequences}. Let $T:\mathcal Z_j\rightarrow
\mathcal Z_j$ again denote the shift.

\mydefinition\label{Definition_J-sequences_Substitutions} For
$y\in X_\sigma$ and $k\in\mathbb Z$, denote by
$[y]_j^k=\{[y(i)]_j\}_{i\in\mathbb Z}$   the
 $j$-sequence  obtained by concatenation
of $j$-symbols $[y(i)]_j$, $i\in\mathbb Z$, where the $j$-symbol
$[y(0)]_j$ appears in $[y]_j^k$ at the position $k$. In other
words, $[y]_j^k$ is a matrix whose $j$-th row is a concatenation
of boxes labeled by $\sigma^j(y(i))$, $i\in\mathbb Z$, from
left to right such that the box labeled by $\sigma^j(y(0))$
starts at the column $k$.

Denote by $\Omega_j$ the closed shift-invariant subset of
$\mathcal Z_j$ generated by all
$j$-sequences  $[y]_j^k$ with $y\in X_\sigma$ and $k\in\mathbb
Z$.\medbreak

The proof of the following lemma is analogous to that of
Proposition \ref{PropositionNumberMinimalComponentsGeneralSubst},
so we omit it.

\begin{lemma}\label{Lemma_Recognaz_NumberMinimalComponents} $(\Omega_j,T)$ has no more than $|A|$ minimal components for any $j\geq 0$.
\end{lemma}

Now we are ready to show that each  substitutional dynamical system without periodic
points is recognizable. Note that the proof uses the ideas of that of
Theorem \ref{TheoremFiniteRank}.

\begin{theorem}\label{TheoremRecognizability} Each aperiodic substitution $\sigma : A\rightarrow A^+$  is recognizable.
\end{theorem}
\proc{Proof.} Assume that $\sigma$ is not recognizable. Then we
have a situation which is shown on the picture below, i.e., there exist $y_{-1},x_{-1}\in X_{\sigma}$, $j_y\in
\{0,\ldots,|\sigma(y_{-1}(0))|-1\}$, and $j_x\in
\{0,\ldots,|\sigma(x_{-1}(0))|-1\}$ with
$y_0=T^{j_y}\sigma(y_{-1})=T^{j_x}\sigma(x_{-1})=x_0$ and  $j_y\neq j_x$ or
$y_{-1}\neq x_{-1}$.

\begin{center}{\small $$\begin{array}{c|c|c|c|c|c|c|c|c|c|c}
\hline \ldots & x(-3) & x(-2) & x(-1) & x(0) & x(1) & x(2) & x(3)
& x(4) & x(5) & \ldots \\
\hline \multicolumn{2}{r|} \ldots & \multicolumn{2}{c|}
{\sigma(y_{-1}(-1))} & \multicolumn{4}{c|} {\sigma(y_{-1}(0))} &
\multicolumn{2}{c|} {\sigma(y_{-1}(1))} & \ldots\\
\hline \ldots & \multicolumn{2}{|c|} {\sigma(x_{-1}(-1))} &
\multicolumn{2}{|c|} {\sigma(x_{-1}(0))} & \multicolumn{4}{|c|}
{\sigma(x_{-1}(1))} & \multicolumn{2}{l} \ldots \\ \hline
\end{array}$$ }
\end{center}

\bigskip (1) Assume that $y_{-1}=x_{-1}$. Hence $j_y\neq j_x$. Thus,
$T^{j_x}\sigma(y_{-1})=T^{j_y}\sigma(y_{-1})$. This shows that
$\sigma(y_{-1})$ is a periodic point, which is impossible.

\medskip (2) Thus, $y_{-1}\neq x_{-1}$. By Proposition \ref{Proposition-1-cutting}, we can find for each
$i\geq 1$ some $x_{-i},y_{-i}\in X_\sigma$ such that
$T^{k_i}\sigma(x_{-i})=x_{-i+1}$ and
$T^{l_i}\sigma(y_{-i})=y_{-i+1}$ for some $l_i\in
\{0,\ldots,|\sigma(y_{-i}(0))|-1\}$ and $k_i\in
\{0,\ldots,|\sigma(x_{-i}(0))|-1\}$. Observe that $x_{-i}\neq
y_{-i}$, $i\geq 1$. Set also $y_i=x_i=\sigma^i(x_0)$, $i\geq 1$.

\medskip (3) Set $K=|A|$ and $L= K^2+K$. It follows from  Proposition
\ref{PropositionSymbolicFormasimNoCommonCuts} that there are  integers
$$0\leq i_0<j_0\leq i_1< j_1\leq\ldots \leq i_{L}<j_{L}=M$$ such
that any pair of $M$-sequences from $\Omega_M$ with depth $i_l$
has no common $j_l$-cuts and, therefore, no common $M$-cuts,
$l=0,\ldots,L$.

For each $i=0,\ldots, M$, define $M$-sequences
$$X_i=[x_{i-M}]_M^{n_i}\mbox{ and }Y_i=[y_{i-M}]_M^{m_i}$$ (see Definition \ref{Definition_J-sequences_Substitutions})
where $n_i$ and $m_i$ are unique integers such that
$$T_\sigma^{n_i}\sigma^{M}(x_{i-M})=x_i=
y_i=T_\sigma^{m_i}\sigma^{M}(y_{i-M}).$$

In other words, the $M$-sequence $X_i$ is built by the rule: the
row 0 consists of concatenated boxes $x_i(j)$, $j\in\mathbb Z$;
 the row 1 consists of concatenated boxes $\sigma(x_{i-1}(j))$
aligned in such a way that if we decompose the row 1 into the letters
we get the equality between rows 0 and 1; and so on; the bottom
line is the concatenation of boxes $\{\sigma^M(x_{i-M}[j])\}$,
$j\in\mathbb Z$, such that whenever we decompose the line $M$ into the
letters we get an equality between rows 0 and $M$. Note that the
$M$-symbol $\sigma^M(x_{i-M}[0])$ appears in $X_i$ at the position
$-n_i$. The $M$-sequence $Y_i$ is built up by the same rule.

\medskip 
(4) It is not hard to see  that for each $l=0,\ldots, L$
the pair $(X_{i_l},Y_{i_l})$  has the depth $i_l$ and, therefore, has no common $j_l$-cuts.

Let $E=E(\Omega_M,T)$ be the enveloping semigroup of
$(\Omega_M,T)$, see Definition
\ref{DefinitionEnvelopingSemigroup}. By definition,
$Ex=\{\gamma(x) : \gamma\in E\}=\overline{Orb_T(x)}$. It follows
that $Ex$ contains at least one minimal $T$-component. It follows
from Lemma \ref{Lemma_Recognaz_NumberMinimalComponents}  that
there is at most $K$ minimal $T$-components, say $Z_0,\ldots,
Z_{K-1}$. Fix any $M$-sequence $Q_i\in Z_i$, $i=0,\ldots, K-1$.

For each pair $(X_{i_l},Y_{i_l})$, $l=0,\ldots,L$, find an element
$\gamma_{i_l}\in E$ such that $\gamma(X_{i_l})=Q_j$ for some
$Q_j$, $j=j(i_l)$.

(4-a) Since  $\gamma_{i_l}$ is the pointwise limit of a sequence
$(T^{n_k})$, the sequences $\gamma_{i_l}(X_{i_l})$ and $\gamma_{i_l}(Y_{i_l})$ remain
$i_l$-compatible.   The fact that the pair $(X_{i_l},Y_{i_l})$  has no
common $j_l$-cuts implies that
$(\gamma_{i_l}(X_{i_l}),\gamma_{i_l}(Y_{i_l}))$ have no common
$j_l$-cuts as well. For otherwise, $T^{n_k}(X_{i_l})$ and
$T^{n_k}(Y_{i_l})$ have a common $j_l$-cut for all $k$ big enough,
which is impossible.

(4-b) Consider the pairs
$(\gamma_{i_l}(X_{i_l}),\gamma_{i_l}(Y_{i_l}))$, $l\in [0,K^2+K]$.
As every $\gamma_{i_l}(X_{i_l})$ coincides with some of
$\{Q_0,\ldots,Q_{K-1}\}$, we can choose $K^2$ pairs
$(\gamma_i(X_{i_l}),\gamma_i(Y_{i_l}))$, $l\in I$, such that
$\gamma_{i_l}(X_{i_l})$ are the same for all $l\in I$, $|I|=K^2$,
say $\gamma_{i_l}(X_{i_l})=Q_{k_0}$ when $l\in I$.

Since each pair $(\gamma_{i_l}(X_{i_l}),\gamma_{i_l}(Y_{i_l}))$ is
$i_l$-compatible and has no common $j_l$-cuts, the choice of the
set $I$ and integers $0\leq i_0<j_0\leq\ldots\ldots i_L<j_L$
guarantees us that each pair
$(\gamma_{i_l}(Y_{i_l}),\gamma_{i_{l'}}(Y_{i_{l'}}))$, $l,l'\in I$ and
$l<l'$, is 0-compatible and has  no common $j_{l'}$-cuts.

(4-c) Setting $Q_l'=\gamma_{i_l}(Y_{i_l})$ for $l\in I$, we obtain a
family $\{Q_l' : l\in I\}\cup\{Q_{k_0}\}$ of  $K^2+1$
$M$-sequences which are $0$-compatible  and have no  common
$M$-cuts. Applying Proposition \ref{PropositionInfectionLemma}, we
get that $\pi_{0}(Q_{k_0})$ is eventually periodic. It follows
from the minimality of $Q_{k_0}$ that $\pi_{0}(Q_{k_0})\in
X_\sigma$ is periodic, which is a contradiction. \ep
\medbreak

For each $a\in A$, set $[a]=\{x\in X_\sigma : x_0=a\}$. Note that,
in general, the set $[a]$ could be
 empty. However, the following result is true even in the case when some of the sets $[a]$, $a\in A$, are empty.

\begin{corollary}\label{CorollaryK-R-partitions}  Let $\sigma : A\rightarrow A^+$ be an aperiodic substitution.
 Then  for every
 $n\geq 0$ $$\mathcal P_n=\{T^i_\sigma\sigma^n([a])\;:\;a\in A\mbox{ and }0\leq i
 <|\sigma^n(a)|\}$$ is a clopen partition of $X_\sigma$.
 Furthermore,  the sequence of partitions
$\{\mathcal P_n\}$ is nested.
 \end{corollary}
\proc{Proof.} First of all, observe that $\sigma^n: X_\sigma
\rightarrow X_\sigma$ is a continuous map.  It follows from
Proposition \ref{Proposition-1-cutting} that the closed set
$\sigma^n(X_\sigma)$ meets each $T_\sigma$-orbit and  consists of
recurrent points. Note that the return time of each $\sigma^n(y)$,
$y\in X_\sigma$ to $\sigma^n(X_\sigma)$ is at most
$|\sigma^n(y[0])|$. Theorem \ref{TheoremRecognizability} implies
that the map $\sigma^n : X_\sigma \rightarrow X_\sigma$ is
one-to-one and the first return time of each $\sigma^n(y)\in
\sigma^n(X_\sigma)$ to $\sigma^n(X_\sigma)$ is exactly
$|\sigma^n(y[0])|$. This shows that $\mathcal P_n$ is a finite
partition of $X_\sigma$ into closed sets and, therefore, into
clopen sets. The fact that the sequence of partitions  $\{\mathcal
P_n\}$ is nested  is proved in \cite[Proposition
14]{durand_host_scau:1999}. \ep

\myremark Observe that, in general, the sequence of K-R partitions
$\{\mathcal P_n\}$ may not generate the topology of $X_\sigma$.
\medbreak

%
%
%
%

%
%
%
%
%
%

\section{Stationary Bratteli-Vershik systems  vs. aperiodic substitutions}\label{SectionDiagramsForSubstitutions}

In the section we show that the class of aperiodic substitutional
systems coincide with the class  of expansive Vershik maps of
stationary ordered Bratteli diagrams.

We recall that Bratteli-Vershik models of substitutional dynamical
systems were constructed  for primitive substitutions in  the
papers \cite{forrest:1997} and \cite{durand_host_scau:1999}. We
also refer the reader to the papers \cite{carlsen_eilers:2004},
\cite{carlsen_eilers:2006}, \cite{yuasa:2002}, and
\cite{yuasa:2006}, where related topics such as various dimension
groups and invariant measures for substitutional dynamical systems
are considered. It is worthwhile to mention the pioneering  paper by
Ferenczi \cite{ferenczi:2005} where the study of substitutions on
infinite alphabets was initiated. We observe that these systems can
be thought as aperiodic homeomorphisms of zero-dimensional Polish
spaces.

\subsection{From Bratteli diagrams to substitutional systems.}
 We start this subsection with the definition of a stationary Bratteli diagram.

\mydefinition\label{Definition_StationaryBratteliDiagram} (1) A
Bratteli diagram $B=(V,E)$ is {\it  stationary} if
$k=|V_1|=|V_2|=\ldots$ and if (by an appropriate labeling of the
vertices) the incidence matrix between levels $n$  and $n+1$ is
the same $k\times k$ matrix $C$ for all $n=1,2,\ldots$. In other
words, beyond level 1 the diagram repeats. Clearly, we can label
the vertices in $V_n$ as $V_n(a_1),\ldots, V_n(a_k)$, where
$A=\{a_1,\ldots,a_k\}$ is a set of $k$ distinct symbols.

(2) A Bratteli diagram $B=(V,E,\leq)$ is {\it stationary ordered}
if $(V,E)$ is stationary and the ordering `$\leq$' on the edges with range
$V_n(a_i)$ does not depend on $n\geq 2$, for all $i=1,\ldots,k$.

(3) Let $B=(V,E,\leq)$ be a stationary ordered Bratteli diagram and   $V_n$ denote the set of
vertices at level $n$, $n\geq 0$. Choose a stationary labeling of
$V_n$  by an alphabet $A$, i.e. $V_n=\{V_n(a)\; |\;
a\in A \}$ for $n>0$.  For every letter $a\in A$, consider the ordered set $(e_1,\ldots,e_k)$ of edges that range at $V_n(a)$, $n\geq 2$, and let
$(a_1,\ldots,a_k)$ be the ordered set of the labels of the
sources of these edges with respect to the ordering `$\leq$'. The map $a \mapsto a_1\ldots a_k$ from $A$
to $A^+$ does not depend on $n$ and therefore determines a
substitution called  the {\it substitution read on $B$}. \medbreak

The following result shows that stationary diagrams can have only a finite number of minimal and maximal paths. See also Propositions \ref{PropositionNumberOfMinimalComponentsVershikMap} and \ref{PropositionNumberMinimalComponentsGeneralSubst}  for similar results.

\begin{proposition}\label{PropositionStationaryOrdering} Let $B=(V,E,\leq)$ be a  stationary ordered Bratteli
diagram.  Then  $B$ has a finite number of minimal and maximal
paths.
\end{proposition}
\proc{Proof.}  Let $\sigma$ be the substitution read on $B$. Set
$$A_0=\{a\in A\; :\; \mbox{ there is }n>0\;(\sigma^n(a)\mbox{ begins with }a)\}$$
and
$$A_1=\{b\in A\; :\; \mbox{ there is }n>0\;(\sigma^n(b)\mbox{ ends with }b)\}.$$
For each $a\in A_0$, let $n_a$ be an integer such that
$\sigma^{n_a}(a)$ begins with $a$. Analogously, for each $b\in
A_1$, let $n_b$ be an integer such that $\sigma^{n_b}(b)$ ends
with $b$.

Setting $p=\prod_{a\in A_0}n_a\times\prod_{b\in A_1}n_b$, we see
that for any $c\in A$ the word $\sigma^n(c)$ begins (ends) with
$c$ for some $n>0$ if and only if $\sigma^p(c)$ begins (ends) with
$c$.

 Now consider a
minimal path $x$ of $B$. Since we have exactly $|A|$
vertices at each level, we can find a letter $a\in A$ and an
infinite set $I\subseteq\mathbb N$ such that $x$ goes through the
vertex $a$ at levels $k$ with $k\in I$. In particular, this means that $\sigma^p(a)$ begins with $a$. We can write
down each $k\in I$ as
$$k=pl_k+m_k\mbox{ with }l_k\geq 0\mbox{ and }0\leq
m_k<p.$$ Find an infinite set $J\subseteq I$ such that
$m_k=m_{k'}$ for $k,{k'}\in J$. It follows that $x$ goes through
the vertex labeled by $a$ at levels $m+np$, $n\geq 1$, where
$m=m_k$ for some $k\in J$. This, in particular, shows that there is
only a finite number of minimal paths.  The proof of the result
for maximal paths is analogous. \ep

\medskip The following proposition shows that  the expansiveness of
Vershik maps defined on stationary Bratteli diagrams is already
seen at the first level.

\begin{proposition}\label{PropositionStationaryCompatible} Suppose that $B=(V,E,\leq)$ is a stationary
Bratteli diagram with  continuous Vershik map. If $x$ and $y$
from $X_B$ are $1$-compatible, then for any $i\geq 1$ there exists
a pair of $i$-compatible elements. In particular, this means that
such a Vershik map is not expansive.
\end{proposition}
\proc{Proof.} We consider a map $f: X_B\rightarrow X_B$, which was
originally defined in \cite{durand_host_scau:1999}. For $x\in
X_B$, let $x_n$ be the label of the edge between levels $n-1$ and
$n$ the path $x$ goes through. For all $n\geq 3$, let $y_n$ be
$x_{n-1}$ and $(y_1,y_2)$ be the minimal path connecting $v_0$ and
$s(y_3)$. Set $f(x)=y$, where $y=(y_1,y_2,\ldots)$.

Now it is not hard to see that if $z_1$ and $z_2$ are
$1$-compatible, then $f^i(z_1)$ and $f^i(z_2)$ are $i$-compatible.
\ep

\medskip The following theorem generalizes the main result of \cite{durand_host_scau:1999} to any aperiodic Vershik map.

\begin{theorem}\label{TheoremVershikMapIsSubstitution} Suppose
that $(X_B,\varphi_B)$ is an aperiodic Bratteli-Vershik system with $B$ a stationary ordered Bratteli diagram and $X_B$ is perfect.
Then the system $(X_B,\varphi_B)$ is homeomorphic to an aperiodic
substitutional dynamical system if and only if no  restriction of
$\varphi_B$ to a minimal component is isomorphic to an odometer.
\end{theorem}
\proc{Proof.} {\it (I)}  Assume first that the diagram $B$ has only
single edges between the top vertex and the vertices of the first level. Let $A$ be the labeling of vertices of $B$ as in Definition \ref{Definition_StationaryBratteliDiagram} and
$\sigma : A\rightarrow A^+$ be the substitution read on $B$.
Define a map $\pi: X_B\rightarrow A^\mathbb Z$ as follows
$$\pi(x)_k=a\mbox{ if and only if }\varphi_B^k(x)
\mbox{ passes through }V_1(a).$$
Observe that $\pi=\pi_1$, where $\pi_1$ maps each $x\in X_B$ to
the row 1 of its matrix, see Section
\ref{SectionExpansiveFiniteRank}. Clearly, $\pi$ is continuous and
\begin{equation}\label{EquationShiftCommuting}\pi\circ
\varphi_B=T\circ \pi,\end{equation} where $T$ is the shift on
$A^\mathbb Z$.

(1)  We claim that $\pi (X_B)\subseteq X_\sigma$. Indeed, to check
this, it is sufficient  to prove that for every $x\in X_B$ and any
$n$ one has $\pi(x){[-n,n]}\in L(\sigma)$.

Assume that $x\in X_B$ is  cofinal neither to a maximal nor to a
minimal path. Let $[x]$ denote the matrix that is obtained by  concatenation of $j$-symbols determined by $x$, see Section
\ref{SectionExpansiveFiniteRank}. Then there is $j>0$ such that
the $j$-symbol, say $v$, from $\mathcal [x]$ crossing 0 column
`covers' coordinates $[-n,n]$ of the first line. This implies that
$\pi(x){[-n,n]}$ is a factor of $\sigma^j(v)$.

By Proposition
\ref{PropositionStationaryOrdering} the diagram has only a finite
number of maximal and minimal paths. Since $X_B$ is perfect, the
set $X_B\setminus Orb_{\varphi_B}(X_{\max}\cup X_{\min})$ is dense
in $X_B$. It follows from the continuity of $\pi$ that
$\pi(X_B)\subseteq X_\sigma$.

(2) We assert  that $\pi(X_B)$ is dense in $X_\sigma$. Indeed,
consider the cylinder set $C=\{x\in X_{\sigma} : x{[n,n+|\sigma^j(a)|-1]}=\sigma^j(a)\}$
with   $a\in A$, $j\geq 1$, and $n\in\mathbb Z$. By Proposition \ref{Proposition-1-cutting}, such cylinder
 sets generate the topology on $X_\sigma$. Thus, we need to show that there is $x\in X_B$ such that $\pi(x)\in C$. Observe that by
(\ref{EquationShiftCommuting}), it is sufficient to find $x\in X_B$
and $k\in \mathbb Z$ such that
$\pi(x){[k,k+|\sigma^j(a)|-1]}=\sigma^j(a)$.

Let $(y_1,\ldots,y_j)$ stand for the minimal path connecting $v_0$
to $V_j(a)$. If $x\in U(y_1,\ldots,y_j)$, then the $j$-symbol
$[a]_j$ appears in the matrix $[x]$   at the coordinate 0.
Therefore, $\pi(x){[0,|\sigma^j(a)|-1]}=\sigma^j(a)$.

(3) It follows from the continuity of $\pi$ that
$\pi(X_B)=X_\sigma$. Thus, $(X_\sigma,T_\sigma)$ is a factor of
$(X_B,\varphi_B)$.

 (4) Suppose that the restriction of $\varphi_B$ onto one of its minimal components is  homeomorphic to an odometer. Then $\varphi_B$ is obviously non-expansive\footnote{It can be deduced, in particular, from Theorem 5.23 of \cite{walters:book}
and the existence of Kakutani-Rokhlin partitions for the system that generate the topology.}. Due
 to the fact that each substitutional dynamical system is
 expansive, the homeomorphism $\varphi_B$ cannot be homeomorphic to
a substitutional systems.

Conversely, if none of the restrictions of $\varphi_B$ to minimal components is homeomorphic to an odometer, then by Theorem \ref{TheoremFiniteRank}
$(X_B,\varphi_B)$ is expansive. By Proposition
\ref{PropositionStationaryCompatible} the map $\pi$ is
injective. This shows that $(X_B,\varphi_B)$ and
$(X_\sigma,T_\sigma)$ are homeomorphic.

\medskip {\it (II)}   Now, let $B$ be an arbitrary ordered Bratteli diagram with the  expansive aperiodic Vershik map $\varphi_B:X_B\rightarrow X_B$. There are two ways of realization of $(X_B,\varphi_B)$ as a substitutional dynamical system. The first one  is to build a  diagram equivalent to $B$  that has only simple edges on the first level, and then apply part (I). The second one is to construct the substitutional dynamical system directly from the diagram $B$. We will exploit the second approach which was first applied in \cite[Proposition 23]{durand_host_scau:1999} for primitive substitutions.

(1) Denote by $A$ the labeling of vertices from $V_n$ ($n\geq 1$) as
in Definition \ref{Definition_StationaryBratteliDiagram}. For each vertex $a\in A$, let $n_a$ be the number of edges
between the vertex $V_1(a)$ of the first level labeled by $a$ and
the top vertex $v_0$.

Denote by $D$ an ordered Bratteli diagram that coincides with $B$ everywhere, but has only one edge
between each vertex $V_1(a)$ and the
top vertex $v_0$. Let $\varphi_D$ be the restriction of
the Vershik map $\varphi_B$ to $X_D$ considered as a clopen subset of $X_B$. Note that $(X_D,\varphi_D)$ is
an induced system of $(X_B,\varphi_B)$. It was proved in (I) that
$(X_D,\varphi_D)$ is homeomorphic to the substitutional dynamical
system $(X_\sigma,T_\sigma)$, where $\sigma$ is the substitution
read on $D$.

(2)   Set $m_a=|\sigma(a)|$, $a\in A$. Since the Vershik map
$\varphi_D$ is aperiodic, we obtain  that $|\sigma^n|\to\infty$ where the substitution $\sigma$ is read on the diagram $D$. Thus, we
can assume, substituting a power of $\sigma$ for $\sigma$ if
needed, that $m_a\geq n_a$ for every $a\in A$.

Let $(X_1,T)$ denote the factor of $(X_B,\varphi_B)$ obtained by
1-sequence coding of elements of $X_B$, see Section
\ref{SectionExpansiveFiniteRank}. Let
$$B_1=\{[a]_1(i) : a\in
A,\;i=0,\ldots,n_a-1\}.$$
One can see that $(X_1,T)$ is a subshift
over the alphabet $B_1$. By Proposition
\ref{PropositionStationaryCompatible}, the factor map $\pi_1 :
X_B\rightarrow X_1$ is injective, so the systems $(X_1,T)$ and
$(X_B,\varphi_B)$ are homeomorphic. Consider the clopen set of $X_1$
$$U=\bigcup_{a\in A}[a](0).$$ Clearly, the return time of any point from
$[a]_1(0)$ to $U$ is $n_a$. We observe  that  $(X_\sigma,T_\sigma)$ is
homeomorphic to $(X_1|_U,T_U)$ and the
homeomorphism is implemented by the map $\xi : A\rightarrow B_1^+$
where $\xi(a)=[a]_1(0)\ldots [a]_1(n_a-1)$ .

(3)  Define the alphabet $E:=\{(a,i) : a\in
A,\; 0\leq i\leq n_a-1\}$ and the map $\psi : A\rightarrow E^+$  by $$\psi(a)=(a,0)(a,1)\ldots
(a,n_a-1).$$
Let $\tau : E\rightarrow E^+$ be a substitution given by

$$\tau(a,i)=\left\{\begin{array}{ll}\psi(\sigma(a)_i) & \mbox{if }1\leq i< n_a-1\\
\psi(\sigma(a)_{[n_a-1,m_a)}) & \mbox{if }i= n_a-1.
\end{array}\right.$$
 Take any $a\in A$,
then $\tau(\psi(a))=\tau(a,0)\ldots
\tau(a,n_a-1)=\psi(\sigma(a))$. Therefore,
\begin{equation}\label{Equation_BratteliDiagrams}\tau^n\circ \psi=\psi\circ \sigma^n\mbox{ for any
}n\geq 0.\end{equation}

For the  substitution $\sigma$ read on the diagram and every
$a\in A$ there are $b\in A$ such that $a$ occurs in $\sigma(b)$,
say, at position $i$. Therefore, $\psi(a)=\psi(\sigma(b)_i)$ is
a factor of $\tau((b,i))$. Thus, $\psi(a)\in L(\tau)$ for every
$a\in A$. By (\ref{Equation_BratteliDiagrams}), we have that
$\psi(L(\sigma))\subseteq L(\tau)$. This shows
that $\psi(X_\sigma)\subset X_\tau$.

Consider the clopen subset $Q$ of $X_\tau$ given by $$Q=\bigsqcup_{a\in
A}[\psi(a)].$$ Clearly, $\psi(X_\sigma)\subseteq Q$.

It follows from Proposition \ref{Proposition-1-cutting} that for every $x\in X_\tau$ there are $z\in X_\tau$ and $0\leq k<|\tau(z[0])|$ such that
$x=T_\tau^k\tau(z)$. Therefore, by definition of $\tau$,
 for every $x\in X_\tau$ there exist
 $y\in A^\mathbb Z$ and  $0\leq k<|\psi(y[0])|$ such
that $x=T^k_\tau\psi(y)$. The definition of $\psi$ implies that such $y$ and $k$ are unique.

Consider any $x\in Q$. Take the unique $y\in A^\mathbb Z$ such that
$\psi(y)=x$. It follows from the definition of $X_\tau$ that every
$\psi(y[-k,k])$ is a factor of $\tau^n((a,i))$ for some $(a,i)\in
E$. Let $a=\sigma(b)_j$ for some $b\in A$ and $j=0,\ldots,m_b-1$, i.e., $(a,i)\prec \psi(\sigma(b))$.
It follows that  $\tau^n((a,i))$ is
a factor of $\tau^{n}\circ \psi(\sigma(b))=\psi(\sigma^{n+1}(b))$.
Therefore, $$\psi(y[-k,k])\prec \psi(\sigma^{n+1}(b)).$$

The definition of $\psi$ implies that if $v\not\prec
w$ where $v,w\in A^+$, then $\psi(v)\not\prec\psi(w)$. Hence $y[-k,k]$ is a factor of $\sigma^{n+1}(b)$. Therefore,  $y\in
X_\sigma$ and $\psi(X_\sigma)=Q$. This shows that $\psi$
implements an homeomorphism between   $(X_\sigma,T_\sigma)$ and
$(X_\tau|_Q,(T_\tau)_Q)$. The following diagram illustrates the
relation between the systems. The downward arrows shows that the
underlying system is an induced system of the overlying one.

 $$\begin{array}{ccccccc} & & (X_1,T) & \stackrel{\pi_1}{\cong} & (X_B,\varphi_B) & \cong & (X_\tau,T_\tau)\\
\multicolumn{2}{c}{\mbox{}}& \downarrow & & \downarrow & &
\downarrow
\\  (X_D,\varphi_D) &  \cong & (X_1|_U,T_U) & \stackrel{\xi}{\cong} &
(X_\sigma,T_\sigma) & \stackrel{\psi}{\cong} &
(X_\tau|_Q,(T_\tau)_Q)\end{array}$$

Observe that the return time of all points from $[\psi(a)]$, $a\in
A$, to $Q$ under the action of $T_\tau$ is $n_a$.    Thus, the map
$\rho : E\rightarrow B_1$ given by $\rho((a,i))=[a]_1(i)$
implements an homeomorphism between $(X_\tau,T_\tau)$ and
$(X_1,T)\cong (X_B,\varphi_B)$. \ep


%

\subsection{From substitutional dynamical systems to Bratteli diagrams.}

In this subsection we show how one can construct a Bratteli-Vershik model for
aperiodic substitutional dynamical systems which satisfy the {\it nesting property} (see Definition \ref{NestingProperty}).
We present a technique applicable for a wide class of substitutions including those with
$|\sigma^n|\to\infty$ and various Chacon-like substitutions.

Let $\sigma :A\rightarrow A^+$ be a substitution. Denote by
$A_l$ the set of all letters $a\in A$ such that
$|\sigma^n(a)|\to\infty$ as $n\to\infty$. Set also $A_s=A\setminus
A_l$.

\mydefinition\label{NestingProperty} We say that a substitution
$\sigma : A\rightarrow A^+$ has a {\it nesting\footnote{We use the
term `nesting' as a natural sequence of K-R partitions for these
substitutional systems is nested, see the proof of Theorem
\ref{TheoremFromSubtitutionsToDiagrams}.} property} if either (1)
for every $a\in A_l$ the word $\sigma(a)$ starts with a letter
from $A_l$ or (2) for every $a\in A_l$ the word $\sigma(a)$ ends
with a letter from $A_l$. \medbreak

\begin{theorem}\label{TheoremFromSubtitutionsToDiagrams}
Let $\sigma : A\rightarrow A^+$ be an aperiodic substitution with
 nesting property. Then the substitutional dynamical system
$(X_\sigma,T_\sigma)$ is homeomorphic to the Vershik map of a
stationary Bratteli diagram.
\end{theorem}
\proc{Proof.} For an alphabet $B$, let $B^*$ denote the set of all
words over $B$ including the empty word.

(1) By Proposition \ref{Proposition_BoundednceOfShortWords}, find
$M>0$ such that any word $W\in L(X_\sigma)$ with $|W|\geq M$
contains a letter from $A_l$.
 Set
$$\mathcal V=\{W\in L(X_\sigma)\; :\; w=v_1S_1v_2S_2v_3,\;v_1,v_2,v_3\in
A_l,\mbox{ and }S_1,S_2\in A_s^*\}.$$ Since no block $S\in
A_s^+\cap L(X_\sigma)$ can have the length greater than $M$, the
set $\mathcal V$ is finite.

Since the substitution $\sigma$ has the nesting property, we have that either (i) all words $\sigma(a)$ start with  letters from $A_l$ for every $a\in A_l$, or (ii) all words $\sigma(a)$ end with  letters from $A_l$ for every $a\in A_l$.

In the first case we define the sequence of K-R partitions as follows
$$\Xi_n=\{T^i_\sigma\sigma^n([v_1S_1.v_2S_2v_3])\; :\; v_1S_1v_2S_2v_3\in\mathcal V,\; 0\leq
i<|\sigma^n(v_2S_2)|\},\; n\geq 0.$$
In the second case we set
$$\Xi_n'=\{T^i_\sigma\sigma^n([v_1.S_1v_2S_2v_3])\; :\; v_1S_1v_2S_2v_3\in\mathcal V,\; 0\leq
i<|\sigma^n(S_1v_2)|\},\;n\geq 0.$$

We will consider in the proof of the theorem the first case only.
The other case  is proved analogously.

(2) We claim that $\{\Xi_n\}$ is a nested sequence of clopen K-R
partitions. Indeed,    assume that
$$
T^i_\sigma\sigma^n([v_0v_1\ldots v_s|v_{s+1}\ldots v_{q-1}])\cap
T^j_\sigma\sigma^n([w_0\ldots w_k|w_{k+1}\ldots w_{r-1}])\neq
\emptyset,
$$
where $0\leq i< |\sigma^n(v_{s+1}\ldots v_{q-2})|$, $0\leq
j<|\sigma^n(w_{k+1}\ldots w_{r-2})|$, the letters  $v_0,v_{s+1},v_{q-1},w_0,w_{k+1},w_{r-1}$ belong to $A_l$, and the
remaining letters $v_t,w_t$ are taken from $A_s$.

Then there exist $y\in [v_0\ldots
v_s|v_{s+1}\ldots v_{q-1}]$ and $z\in [w_0\ldots w_k|w_{k+1}\ldots
w_{r-1}]$ with $T^i\sigma^n(y)=T^j\sigma^n(z)$. By Theorem
\ref{TheoremRecognizability} and the definition of $\mathcal V$,
we get that $y=z$ and $i=j$. This implies that $q=r$, $k=s$, and
$v_t=w_t$ for $t=0,\ldots,q-1$. Thus, $\Xi_n$ is a  K-R partition of
$X_\sigma$. Notice that, by  continuity of $\sigma^n$, all the sets
$T^i_\sigma\sigma^n([v_1S_1.v_2S_2v_3])$,
$v_1S_1v_2S_2v_3\in\mathcal V$, are closed.  Therefore, they are clopen.

 Applying $\sigma$ to any
$W=v_1S_1|v_2S_2v_3\in\mathcal V$, we get a representation
$$
\sigma(W)=w_1^{(1)}F_1^{(1)}\ldots w_{k_1}^{(1)}F_{k_1}^{(1)}|w_{1}^{(2)}F_{1}^{(2)}\ldots w_{k_2}^{(2)}F_{k_2}^{(2)} w_{1}^{(3)}F_{1}^{(3)}\ldots w_{k_3}^{(3)}F_{k_3}^{(3)}
$$ with $w_i^{(\cdot)}\in A_l$, and $F_i^{(\cdot)}\in A_s^*$. Here $\sigma(v_iS_i)=w_{1}^{(i)}F_{1}^{(i)}\ldots w_{k_i}^{(i)}F_{k_i}^{(i)}$ where $i=1,2,3$ and $S_3=\emptyset$. Since $\sigma$ has a nesting property, we obtain that $w_{1}^{(i)}\neq\emptyset$. This representation
shows that the $T_\sigma$-tower from $\Xi_n$ with base
$\sigma^n([v_1S_1.v_2S_2v_3])$ consequently intersects the
$T_\sigma$-towers of $\Xi_{n-1}$ with bases:
$$
\begin{array}{ll}\sigma^{n-1}([w_{k_1}^{(1)}F_{k_1}^{(1)}.w_{1}^{(2)}
F_{1}^{(2)}w_{2}^{(2)}]),\;
\sigma^{n-1}([w_{2}^{(2)}F_{2}^{(2)}.
w_{3}^{(2)}F_{3}^{(2)}w_{4}^{(2)}]),...,\\
\\
\sigma^{n-1}([w_{k_2-1}^{(2)}F_{k_2-1}^{(2)}.
w_{k_2}^{(2)}F_{k_2}^{(2)}w_{1}^{(3)}]).\end{array}
$$
In particular, we obtain that $B(\Xi_n)\subseteq B(\Xi_{n-1})$ and
$\{\Xi_n\}_{n\geq 0}$ is a nested sequence of K-R partitions.

 (3) We claim that the partitions $\{\Xi_n\}_{n\geq 0}$
generate the topology of $X_\sigma$. Observe that it suffices  to
show that the function $x_{[-m,m]}$, $m>0$, is constant on
each element of partition $\{\Xi_n\}$ for $n$  big enough.
To see this, choose $n$ such that $$\min\{|\sigma^{n}(a)|\; : \;a\in
A_l \}>m.$$

Fix a word $W=(v_1S_1v_2S_2v_3)\in\mathcal V$ and $0\leq
k<|\sigma^n(v_2S_2)|$. For each $x\in
T^k\sigma^n([v_1S_1.v_2S_2v_3])$, there exists $y\in
[v_1S_1.v_2S_2v_3]$ such that  $x=T^k\sigma^n(y)$. It follows that
the word $\sigma^n(v_2S_2v_3)$ is a prefix of
$\sigma^n(y{[0,\infty)})$ and $\sigma^n(v_1S_1)$ is a suffix of
$\sigma^n(y{(-\infty,-1]})$. Therefore,
$$\sigma^n(y){[-L,R)}=\sigma^{n}(v_1S_1v_2S_2v_3),$$
where $L=|\sigma^{n}(v_1S_1)|$ and $R=|\sigma^n(v_2S_2v_3)|$.
Thus, we obtain
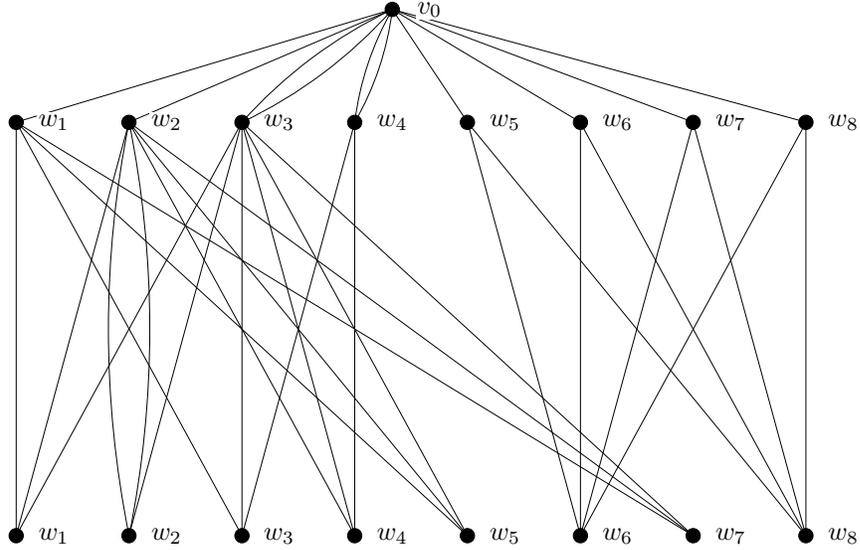
$$x_{[-m,m]}=\sigma^{n}(v_1S_1v_2S_2v_3)_{[L+k-m,L+k+m)},$$
which does not depend on $x$, but only on $k$ and the word $W$.

Now it follows from Theorem \ref{ExistenveBratteliDiagrams}  that
the substitutional system $(X_\sigma,T_\sigma)$ is homeomorphic to
the Bratteli-Vershik system $(X_B,\varphi_B)$, where the ordered
Bratteli diagram $B$ is constructed by the sequence of K-R
partitions $\{\Xi_n\}_{n\geq 0}$. The fact that the diagram $B$ is
stationary  has been proved in (2). \ep

\myexample\label{ExampleExtendedChacon} Consider the Chacon  type substitution on the alphabet $A=\{0,s,1\}$ given by
$\sigma(0)=00s0$, $\sigma(s)=s$, and $\sigma(1)=0110$. Define the
set $\mathcal V$ as in the proof of Theorem
\ref{TheoremFromSubtitutionsToDiagrams}. Then $\mathcal V$ consists of the words $w_1,\ldots,w_8$, where
$$\begin{array}{llll} w_1=0.00 & w_2=0s.00 & w_3=0.0s0 & w_4=0s.0s0\\
w_5=1.00 & w_6=0.11 & w_7=0.01 & w_8=1.10\end{array}$$

Here the dot separates the negative and non-negative coordinates as in
the definition of K-R partitions $\{\Xi_n\}$ in Theorem
\ref{TheoremFromSubtitutionsToDiagrams}. To construct the Bratteli
diagram, we need to trace the orbits of each base
$\sigma^n[w_i]$, $i=1,\ldots,8$. In other words, we need to list
all $T_\sigma$-towers which are intersected by $\sigma^n[w_i]$.
This is achieved by decomposing each word $\sigma(w_i)$,
$i=1,\ldots,8$, into the letters and analyzing the occurrence of long
and short letters. We  clarify the above scheme by considering
the set $\sigma^n[w_4]$. Decompose $\sigma(w_4)$ into the letters:
$$\begin{array}
{c@{\,}c@{\,}c@{\,}c@{\,}c@{\,}c@{\,}c@{\,}c@{\,}c@{\,}c@{\,}c@{\,}
c@{\,}c@{\,}c@{\,}c@{\,}}\cline{8-12}
& & & & & & & \multicolumn{5}{|c|} {w_4}%
\\%
\cline{5-8}
& & & & \multicolumn{3}{|c|}{w_2} & \multicolumn{1}{c|}{} & &
\multicolumn{1}{c}{}& 
& &
\multicolumn{1}{|c}{} &
\\
\cline{7-10}
\sigma(w_4)= &\multicolumn{1}{c}{ 0 }&\multicolumn{1}{c}{ 0 }& \multicolumn{1}{c}{ s }&\multicolumn{1}{|c} 0 &{ s.} &
\multicolumn{1}{|c} 0
&\multicolumn{1}{|c|} 0 & s &\multicolumn{1}{c|} 0 & s & \multicolumn{1}{c}0 & \multicolumn{1}{|c}0 & s & \multicolumn{1}{c}0%
\\
\cline{5-12} & & & & & & \multicolumn{1}{|c} {} & &
\multicolumn{1}{c}{ w_3}& & \multicolumn{1}{|c} {} & & &
\multicolumn{1}{c}{}
\\ \cline{7-10} & & & & & & & & & 
\end{array}$$

This shows that the base $\sigma^n[w_4]$ consequently meets the
bases $\sigma^{n-1}[w_2]$, $\sigma^{n-1}[w_3]$, and 
$\sigma^{n-1}[w_4]$.
We  denote this correspondence by $w_4\mapsto
w_2w_3w_4$. Repeating this argument for every set
$\sigma^n([w_i])$, we get the following matching rule
$$\tau:\left\{\begin{array}{l} w_1\mapsto w_1w_3w_2\\
w_2\mapsto w_2w_3w_2\\
w_3\mapsto w_1w_3w_4\\
w_4\mapsto w_2w_3w_4\\
w_5\mapsto w_1w_3w_2\\
w_6\mapsto w_7w_6w_8w_5\\
w_7\mapsto w_1w_3w_2\\
w_8\mapsto w_7w_6w_8w_5
\end{array}\right.$$

 It follows from the proof of Theorem \ref{TheoremFromSubtitutionsToDiagrams}
 that the system $(X_\sigma,T_\sigma)$ is conjugate to the Vershik map
of  the following stationary ordered Bratteli diagram determined
by the rule $\tau$, i.e., $\tau$ is the substitution read from
$B$.

\unitlength = 0.5cm
\begin{center}
\begin{graph}(26,17)
\graphnodesize{0.4} \roundnode{V0}(12,16)
\nodetext{V0}(1,0){$v_0$}
\roundnode{W11}(2,13)
\roundnode{W21}(5,13)
\roundnode{W31}(8,13)
\roundnode{W41}(11,13)
\roundnode{W51}(14,13)
\roundnode{W61}(17,13)
\roundnode{W71}(20,13)
\roundnode{W81}(23,13)
\roundnode{W12}(2,2)
\roundnode{W22}(5,2)
\roundnode{W32}(8,2)
\roundnode{W42}(11,2)
\roundnode{W52}(14,2)
\roundnode{W62}(17,2)
\roundnode{W72}(20,2)
\roundnode{W82}(23,2)
\edge{V0}{W11}
\edge{V0}{W21}
\bow{V0}{W31}{0.05}
\bow{V0}{W31}{-0.05}
\bow{V0}{W41}{0.05}
\bow{V0}{W41}{-0.05}
\edge{V0}{W51}
\edge{V0}{W61}
\edge{V0}{W71}
\edge{V0}{W81}
\nodetext{W11}(1,0){$w_1$} \nodetext{W21}(1,0){$w_2$}
\nodetext{W31}(1,0){$w_3$} \nodetext{W41}(1,0){$w_4$}
\nodetext{W51}(1,0){$w_5$} \nodetext{W61}(1,0){$w_6$}
\nodetext{W71}(1,0){$w_7$} \nodetext{W81}(1,0){$w_8$}
\nodetext{W12}(1,0){$w_1$} \nodetext{W22}(1,0){$w_2$}
\nodetext{W32}(1,0){$w_3$} \nodetext{W42}(1,0){$w_4$}
\nodetext{W52}(1,0){$w_5$} \nodetext{W62}(1,0){$w_6$}
\nodetext{W72}(1,0){$w_7$} \nodetext{W82}(1,0){$w_8$}

%
\edge{W12}{W11}
\edge{W12}{W21}
\edge{W12}{W31}
%
\bow{W22}{W21}{-0.05}
\bow{W22}{W21}{0.05}
\edge{W22}{W31}
%
\edge{W32}{W11} 
\edge{W32}{W31}
\edge{W32}{W41}
%
\edge{W42}{W21}
\edge{W42}{W31}
\edge{W42}{W41}
%
\edge{W52}{W11}
\edge{W52}{W21}
\edge{W52}{W31}
%
\edge{W62}{W51}
\edge{W62}{W61}
\edge{W62}{W71}
\edge{W62}{W81}
%
%
\edge{W72}{W11}
\edge{W72}{W21}
\edge{W72}{W31}
%
\edge{W82}{W71}
\edge{W82}{W61}
\edge{W82}{W81}
\edge{W82}{W51}
\end{graph}
\end{center}
Notice that by Theorem \ref{TheoremVershikMapIsSubstitution}, the
systems $(X_\sigma,T_\sigma)$ and $(X_\tau,T_\tau)$ are conjugate.
See also Example \ref{ExampleDerivativeSubtition} for another
Bratteli-Vershik model of $(X_\sigma,T_\sigma)$. We  also mention the work \cite[Section 4.2]{gjerde_johansen:2000} where the authors presented a Bratteli-Vershik model for the minimal component of the system $(X_\sigma,T_\sigma)$. \medbreak

%
%
%

\section{Derivative substitutions}\label{Section_DerivativeSubstitutions}
In the section, we study a generalization of the notion of derivative substitutions which was defined
in \cite{durand:1998} and \cite{durand_host_scau:1999} for primitive substitutions. We show  that this notion works also for  non-primitive substitutions.

To make the exposition clear and abandon some pathological
situations,  we  restrict our study to the class of substitutions,
which we call {\it $m$-primitive}. However, the results can be
applied to more general class of substitutions, for example, for
the Chacon type substitutions.

\mydefinition Let $m>0$. We say that a substitution $\sigma:
A\rightarrow A^+$ is {\it $m$-primitive} if  we can decompose the
alphabet $A=A_1\sqcup\ldots \sqcup A_m\sqcup A_0$ with $|A_i|\geq
2$, $i=1,\ldots,m$, and

(a) $\sigma(A_i)\subset A^+_i$ and $\sigma|_{A_i}$ is primitive
for every $i=1,\ldots,m$;

(b) for all $a\in A$ there exists a letter $b\in A_1\cup\ldots\cup
A_m$ and $i>0$ with $b\prec\sigma^i(a)$.

(c) $L(\sigma)=L(\sigma^k)$ for every $k\geq 1$.

(d)  $A\subset L(X_\sigma)$.\medbreak

Given an $m$-primitive  substitution $\sigma :A\rightarrow A^+$,
define the matrix of the substitution
$M(\sigma)=(M(\sigma)_{a,b})_{a,b\in A}$ as follows:
$M(\sigma)_{a,b}$  is the number of occurrences of $b$ in
$\sigma(a)$. Then the matrix $M(\sigma)$ is of the form \begin{equation}M(\sigma)=\left(%
\begin{array}{ccccccccc}
 M_1 & O & O & \ldots & O & O &  \ldots & O \\
  O & M_2 & O & \ldots& O & O & \ldots & O \\
  \vdots & \vdots &\vdots &\ddots &\vdots &\vdots  &\ddots &\vdots \\
  O & O & O & \ldots & M_m & O &  \ldots & O\\
  p_{1,1}   & p_{1,2} & p_{1,3} & \ldots & p_{1,m} & p_{1,m+1} &  \ldots & p_{1,m+q}\\
  p_{2,1}   & p_{2,2} & p_{2,3} & \ldots & p_{2,m} & p_{2,m+1}  & \ldots & p_{2,m+q}\\
  \vdots & \vdots &\vdots &\ddots &\vdots &\vdots &\ddots &\vdots \\
  p_{n,1}   & p_{n,2} & p_{n,3} & \ldots & p_{n,m} & p_{n,m+1}  & \ldots &   p_{n,m+q}
\end{array}%
\right)\end{equation} where $O$ denotes a zero matrix,  the matrices
$M_i$, $i=1,\ldots,m$,  are primitive, and for each $j=1,\ldots,n$ there is a power $k$ of $M$ such that
at least one of the entries $p_{j,1}^{(k)},\ldots, p_{j,m}^{(k)}$ from $M^k$  is not zero. We notice also that $|\sigma^n|\to\infty$.

\medskip  By Proposition
\ref{PropositionNumberMinimalComponentsGeneralSubst} the number of
minimal components of the substitutional system
$(X_\sigma,T_\sigma)$ is bounded by $|A|$. However, for the class
of $m$-primitive substitutions, the set of all minimal components
admits a complete description.

\begin{proposition}\label{PropositionNumberMinimalComponents-M-primitive} Let $\sigma$ be an $m$-primitive substitution.
Then $(X_\sigma, T_\sigma)$ has exactly $m$ minimal components,
which are $$Z_i=\{x\in X_\sigma\; :\; x_{[-n,n]}\mbox{ is a factor
of some }\sigma^m(a),\,a\in A_i\},$$ $i=1,\ldots,m$.
\end{proposition}
\proc{Proof.} The fact that $Z_i$, $i=1,\ldots,m$, are minimal
components is well-known  \cite{queffelec:book}. We will  show
that there are no others. For each $Z_i$,
we can find $w_i\in Z_i$  and $p_i>0$  such
that $\sigma^{p_i}(w_i)=w_i$, see for example \cite[Chapter
V]{queffelec:book}. Without loss of generality, we may assume that
$\sigma(w_i)=w_i$ for all $i=1,\ldots,m$. We  also assume that
all letters from $A_i$ appear in $\sigma(a)$ for every $a\in A_i$,
$i=1,\ldots,m$.

Take any $n>0$. For every $x\in X_\sigma$, there is $a\in A$ such
that the word $\sigma^{n+2}(a)$ appears in $x$.  Find  $b\in A_i$  that appears  in $\sigma(a)$ for some $i=1,\ldots,m$.  Therefore, the
word $\sigma^n(w_i(0))$ appears in $x$.  Find $1\leq i_0\leq m$ and an infinite set $I\subseteq \mathbb N$ such that the
word $\sigma^n(w_{i_0}(0))$ appears in $x$ for all $n\in I$.  This implies
that the closure of the orbit of $x$ contains the point $w_{i_0}$.
Hence, $\overline{Orb_{T_\sigma}(x)}\supseteq Z_{i_0}$. This proves
the result. \ep

\myremark Note that $(X_\sigma,T_\sigma)$ has a periodic point iff
$(Z_i,T_\sigma)$ is periodic for some $i$. Observe also that there
is an algorithm that decides whether $\sigma$ has periodic points
or not, see, for example,  the book \cite{kurka:book} or
references in \cite{durand_host_scau:1999}.

\subsection{Return words and proper substitutions}
Let $\sigma: A\rightarrow A^+$ be an $m$-primitive substitution.
Here we assume that $(X_\sigma,T_\sigma)$ has no periodic points.
By Proposition
\ref{PropositionNumberMinimalComponents-M-primitive}, we get that
the system $(X_\sigma,T_\sigma)$ has exactly $m$ minimal
components $Z_1,\ldots,Z_m$. By standard arguments, we can find
$w_i\in Z_i$  and $p_i>0$  such that $\sigma^{p_i}(w_i)=w_i$  \cite[Chapter V]{queffelec:book}. Without loss of
generality, we can assume that $\sigma(w_i)=w_i$ for all
$i=1,\ldots,m$. Notice also that
$\overline{Orb_{T_\sigma}(w_i)}=Z_i$ for each $i=1,\ldots,m$.

Set $r_i=w_i[-1]$ and $l_i=w_i[0]$ for $i=1,\ldots,m$. Observe
that each word $\sigma(r_i)$ ends with $r_i$, whereas every word
$\sigma(l_i)$ begins with $l_i$, for all $i=1,\ldots,m$.

Consider the clopen sets $[r_i.l_i]:=\{x\in X_\sigma :
x[-1]=r_i,\;x[0]=l_i\}$ and define $$W=\bigcup_{i=1}^m[r_i.l_i].$$

\begin{proposition}\label{PropositionReturnTime}  The set $W$ meets each $T_\sigma$-orbit and consists of recurrent
points. Furthermore, the return time to $W$ is bounded.
\end{proposition}
\proc{Proof.} Take any $x\in X_\sigma$ and consider
$Z(x)=\overline{Orb_{T_\sigma}(x)}$. Since $Z(x)$ is a
$T_\sigma$-invariant closed subset of $X_\sigma$, it contains one of the minimal components $Z_i$, $i=1,\ldots,m$, say $Z_{i_0}$.
Therefore, $[r_{i_0}.l_{i_0}]$ meets the orbit of $x$. Hence
$X_\sigma=\bigcup_{n\in\mathbb Z}T^n_\sigma W$. By  compactness
of $X_\sigma$, we see that $X_\sigma=\bigcup_{n=0}^kT^n_\sigma W$
for some $k>0$. This shows that $W$ consists of recurrent points.
\ep

\mydefinition\label{DefinitionReturnWords} We say that $w\in L(X_\sigma)$ is a {\it return
word} if there exists $1\leq i,j\leq m$ such that

 (i) $r_iwl_j\in L(X_\sigma)$;

(ii) the first and  last letters of $w$ are $l_i$ and $r_j$,
respectively;

(iii) no  word from $\{r_1l_1,\ldots,r_ml_m\}$   appears in $w$.
\medbreak

Let $\mathcal R$ denote the set of all return words. In view of
Proposition \ref{PropositionReturnTime}, the set $\mathcal R$ is
finite.  Let us enumerate the return words in an arbitrary way and
denote $\mathcal N=\{1,2,\ldots,card(\mathcal R)\}$. Let $\phi: \mathcal N\rightarrow \mathcal R$ be an
``enumeration'' map.

The proof of the following proposition is trivial, so we omit it.

\begin{proposition}\label{PropositionDerivSubst1-1} The
maps $\phi_1: \mathcal N^+\rightarrow A^+$ and $\phi_2: \mathcal
N^\mathbb Z\rightarrow A^\mathbb Z$ induced by $\phi$ are one-to-one.
\end{proposition}

Take any return word $w\in\mathcal R$. Decompose it into letters
$w=w_1\ldots w_k$. By definition of return words, $w_1=l_i$ and
$w_k=r_j$ for some $i,j$, and $r_iwl_j\in L(\sigma)$. Therefore,
$\sigma(r_i)\sigma(w)\sigma(l_j)\in L(\sigma)$. Since
$\sigma(r_i)$ ends with $r_i$, whereas every word $\sigma(l_i)$
begins with $l_i$, we have that $r_i\sigma(w)l_j\in L(\sigma)$.
So, the word $\sigma(w)$ appears in the word $r_i\sigma(w)l_j$
between occurrences of $r_i.l_i$ and $r_j.l_j$. Thus, by Proposition \ref{PropositionDerivSubst1-1}, $\sigma(w)=d_1,\ldots,d_q$ can be uniquely written as a concatenation of return words $d_i$.

\mydefinition\label{DerivativeSubstitution} Define the substitution $\tau : \mathcal
R\rightarrow \mathcal R^+$ by $\tau(w)=d_1\ldots d_q$, where
$d_1\ldots d_q$ is the unique decomposition of $\sigma(w)$ into
return words. The substitution $\tau$ is called the {\it derivative
substitution} of $\sigma$.  \medbreak

Notice that
 $$\phi\circ \tau^m=\sigma^m\circ \phi.$$

The following result justifies the name of the derivative
substitution.

\begin{proposition}\label{PropositionTopolInterprDerivSubst} The substitutional dynamical system $(X_\tau,T_\tau)$
associated to $\tau:\mathcal R\rightarrow \mathcal R^+$ is
homeomorphic to the system induced by $(X_\sigma,T_\sigma)$ on the
cylinder set $W=\bigcup_{i=1}^m[r_i.l_i]$.
\end{proposition}
\proc{Proof.} By Proposition \ref{PropositionDerivSubst1-1}, we
have that $\phi:\mathcal R^\mathbb Z\rightarrow A^\mathbb Z $ is
one-to-one.

(1) $\phi(X_\tau)\subset X_\sigma$. Indeed, consider any
$x=(x_i)\in X_\tau$. By definition of $X_\tau$, for every $n>0$,
there are $w\in \mathcal R$ and an integer $m>0$ such that
$x_{[-n,n]}$ is a factor of $\tau^m(w)$. It follows that
$\phi(x_{[-n,n]})$ is a factor of $\phi\circ
\tau^m(w)=\sigma^m\circ \phi (w)$. Since $\phi(w)\in L(\sigma)$,
we get that $\phi(x_{[-n,n]})\in L(\sigma)$ for every $n>0$.

(2)  It is clear that $\phi(X_\tau)\subseteq W$.

(3) $\phi(X_\tau)= W$. Indeed, take any $z\in W$. By Proposition
\ref{PropositionReturnTime}, we can decompose $z$ into a
concatenation of return words $z=\ldots \phi(x[-1]) |\phi(
x[0])\ldots$. We must show that $x=(x[i])\in X_\tau$. By the definition of $X_\sigma$ for any $n>0$ there are
$k>0$ and $a\in A$ such that $\phi(x[-n,n])\prec \sigma^k(a)$. Take a return word $w\in \mathcal R$ such that  the letter $a$ is
a factor of $\phi(w)$.
Therefore, $$\phi(x[-n,n])\prec \sigma^k(a)\prec \sigma^k(\phi(w))=\phi(\tau^k(w)).$$ It follows from the definition of $\phi$ that
$x[-n,n]$ is a factor of $\tau^k(w)$.

(4) If  $x\in X_{\tau}$ and $z=\phi(x)$, then the first return
time of $z$ to $W$ is $n=|\phi(x_0)|$. Thus, the image of $z$ by
the first return time transformation is
$T_{\sigma}^nz=\phi(T_{\tau}x)$. This proves the proposition. \ep

\medskip Now we introduce the notion of a proper substitution. Our
motivation to bring in this notion comes from the paper
\cite{durand_host_scau:1999}, where proper substitutions were
indispensable for the Bratteli diagram construction. However, our
definition of proper substitutions differs from that of
\cite{durand_host_scau:1999}. The key concept which we put
behind this notion is that for proper substitutions (for both
definitions: ours and from \cite{durand_host_scau:1999}) the
sequence of K-R partitions $\{\mathcal P_n\}$ defined in Corollary
\ref{CorollaryK-R-partitions} generates the topology.

\mydefinition\label{Definition_ProperSubstitution} Let $\sigma :
A\rightarrow A^+$ be a substitution. We say that $\sigma$ is {\it
proper} if there is $p>0$ such that for every letter $a\in A$, the
first letter of $\sigma^p(b)$, with $b\in A$ and $ab\in
L(X_\sigma)$, does not depend on $b$; and the last letter of
$\sigma^p(c)$, with $c\in A$ and $ca\in L(X_\sigma)$, does not
depend on $c$. \ep

\myremark Observe that if a substitution is proper in the sense of
\cite{durand_host_scau:1999}, then it is also proper by our
definition. \medbreak

\begin{proposition}\label{PropositionDerivativeSubstitutionIsProper}
The derivative substitution  $\tau$, defined by an $m$-primitive substitution,  is a proper aperiodic $m$-primitive substitution.
\end{proposition}
\proc{Proof.}  Denote by $\mathcal R_i$ the set of all return
words that appear in the fixed point $w_i$, $i=1,\ldots,m$. We will identify the elements of $\mathcal R$ with their counterparts in $L(X_\sigma)$ (see Definition \ref{DefinitionReturnWords}).
 Set also $\mathcal R_0=\mathcal R\setminus (\mathcal R_1\cup\ldots
\cup\mathcal R_m)$. To verify that $\tau :
\mathcal R_i\rightarrow\mathcal R_i^+$ is a primitive substitution,
we refer the reader to the proof of  Lemma 21 from
\cite{durand_host_scau:1999}.

 Denote by $v_i$  the return word that appears first (if we are
counting rightwards) in the fixed point $w_i$, $i=1,\ldots,m$.
That is $v_i$ is a prefix of $(w_i)_{[0,+\infty)}$. Find $n>0$
such that
\begin{equation}\label{Equation_DerivativeISProper}|\sigma^n|>
\max_{i=1,\ldots,m}|v_i|.\end{equation} It
follows that $v_il_i$ is a prefix of $\sigma^n(l_i)$. Take any
return words $w$ and $w'$ with $ww'\in L(\tau)$. By definition,
$r_{i_0}$ is a suffix of $w$ for some $i_0=1,\ldots,m$. Therefore,
the word $w'$ begins with $l_{i_0}$. So $\sigma^n(l_{i_0})$ is a
prefix of $\sigma^n(w')$. This implies that $v_{i_0}l_{i_0}$ is a
prefix of $\sigma^n(w')$. This means that $v_{i_0}$ is the first
return word in $\sigma^n(w')$. That is $\tau^n(w')$ begins with
$v_{i_0}$. Thus, the first letter of $\tau^n(w')$ does not
depend on $w'$, but only on $w$.

The same argument works to show that the last letter of
$\tau^n(w')$ with $w'w\in L(\tau)$ depends only on $w$ for $n$
large enough.

Now consider any $w\in\mathcal R_0$. Then the first letter of $w$
is $l_i$ for some $i=1,\ldots,m$. It follows that $\tau^n(w)$
contains $v_i$, where $n$ is as in
(\ref{Equation_DerivativeISProper}). \ep

\medskip To prove the following result, we use  Proposition 14 from \cite{durand_host_scau:1999}.

\begin{proposition}\label{PropositionProperSubstGeneratingTopology} Let $\sigma : A\rightarrow A^+$ be
an aperiodic proper substitution such that $|\sigma^n|\to\infty$. Then  the sequence of partitions
$\{\mathcal P_n\}$ defined in Corollary
\ref{CorollaryK-R-partitions} generates the topology of $X_\sigma$.
\end{proposition}
\proc{Proof.} Let $p>0$ be an integer as in Definition
\ref{Definition_ProperSubstitution} of proper substitutions. Given
an integer $m>0$, we claim that for $n$ sufficiently large  the function
$x_{[-m,m]}$ is constant on each element of partition $\{\mathcal
P_n\}$. Choose $n$ so large that $|\sigma^{n-p}|>m$.

Fix $a\in A$ and $0\leq k<|\sigma^n(a)|$. For each $x\in
T^k\sigma^n([a])$, there exists $y\in X_\sigma$ such that $y_0=a$
and $x=T^k\sigma^n(y)$. The word $\sigma^n(a)\sigma^n(y_1)$ is a
prefix of $\sigma^n(y_{[0,\infty)})$. By definition of proper
substitutions, the first letter of $\sigma^p(y_1)$, say $l$, does
not depend on $y_1$ (it depends only on the letter $a$). It
follows that $\sigma^n(a)\sigma^{n-p}(l)$ is a prefix of
$\sigma^n(y_{[0,\infty)})$. Similarly, there is $r\in A$
that depends only on the letter $a$ (not on $y$) and such that
$\sigma^{n-p}(r)$ is a suffix of $\sigma^n(y_{(-\infty,-1]})$.
Therefore,
$$\sigma^n(y)_{[-R,L)}=\sigma^{n-p}(r)\sigma^n(a)\sigma^{n-p}(l),$$
where $R=|\sigma^{n-p}(r)|$ and
$L=|\sigma^n(a)|+|\sigma^{n-p}(l)|$. Thus, we obtain
$$x_{[-m,m]}=\sigma^{n-p}(r)\sigma^n(a)\sigma^{n-p}(l)_{[R+k-m,R+k+m)},$$
which does not depend on $x$, but only on $k$ and $a$. \ep

\medskip Now, we are ready to present a general approach for construction of
Bratteli diagrams for $m$-primitive substitutions. To build a
Bratteli-Vershik model for $\sigma$, it is sufficient to construct a
Bratteli diagram for $\tau$ and then add some edges to the first level.

\begin{theorem}\label{TheoremBratteliDiagrForSubstitutions_ViaDerivative} Let $\sigma$ be an $m$-primitive aperiodic substitution over an
alphabet $A$ with  derivative substitution  $\tau$. Let $B=(V,E,\leq)$ be the stationary ordered Bratteli diagram built by the matrix of the substitution $\tau$ where $\tau$ is  read on $B$. Suppose also that $B$
has $|w|$ edges between the top vertex and the vertex defined by
the return word $w\in\mathcal
 R$. Then $B$ admits a continuous dynamics $(X_B,T_B)$ which is  homeomorphic  to $(X_\sigma,T_\sigma)$.
\end{theorem}
\proc{Proof.}  Let $B'$ be the stationary Bratteli diagram built
by the matrix of the substitution $\tau$ with simple edges between the top vertex
and vertices of the first level. By Proposition
\ref{PropositionDerivativeSubstitutionIsProper} $\tau$ is a proper
substitution. It follows from Proposition
\ref{PropositionProperSubstGeneratingTopology} and results of
Section \ref{Section_NotionBratteliDiagram} that the Bratteli
diagram $B'$ admits continuous dynamics $(X_{B'},\varphi_{B'})$
which is homeomorphic to $(X_\tau,T_\tau)$. Then the application of
Proposition \ref{PropositionTopolInterprDerivSubst} yields the
result. \ep

\myexample\label{ExampleDerivativeSubtition} Let $\sigma$ be the
 Chacon type substitution defined in  Example
\ref{ExampleExtendedChacon}. Note that the system
$(X_\sigma,T_\sigma)$ has only one minimal component which is
generated by the fixed point
$$w=\lim_n\sigma^n(0).\sigma^n(0)=\ldots  00s0 s 00s0. 00s0 00s0 s 00s0
\ldots$$ Though $\sigma$ is not an $m$-primitive substitution, the technique developed in the section can still be applied to build a Bratteli-Vershik model of $(X_\sigma,T_\sigma)$.

 Consider the set of all return words
$\mathcal R$. One can check that $$\mathcal R
=\{v_1=0,\;v_2=0s0,\;v_3=0s0s0,\;v_4=0110\}.$$ Find the unique
decomposition of each $\sigma(v_i)$ into return words
$$\begin{array}{l}\sigma(v_1)=0|0s0=v_1v_2\\
\sigma(v_2)=0|0s0 s 0|0s0=v_1v_3v_2\\
\sigma(v_3)= 0|0s0 s 0|0s0 s 0|0s0=v_1v_3v_3v_2\\
\sigma(v_4) = 0|0s0| 0110 |0110| 0|0s0=v_1v_2v_4v_4v_1v_2
\end{array}$$
Thus, the derivative substitution $\tau:\mathcal R\rightarrow
\mathcal R^+$ is defined by $\tau(v_1)=v_1v_2$;
$\tau(v_2)=v_1v_3v_2$; $\tau(v_3)=v_1v_3v_3v_2$; and
$\tau(v_4)=v_1v_2v_4v_4v_1v_2$. Clearly, the substitution $\tau$ is proper and $|\tau^n|\to\infty$. By Theorem
\ref{TheoremBratteliDiagrForSubstitutions_ViaDerivative}, the
system $(X_\sigma,T_\sigma)$ is conjugate to the Vershik map of
the following stationary ordered Bratteli diagram:

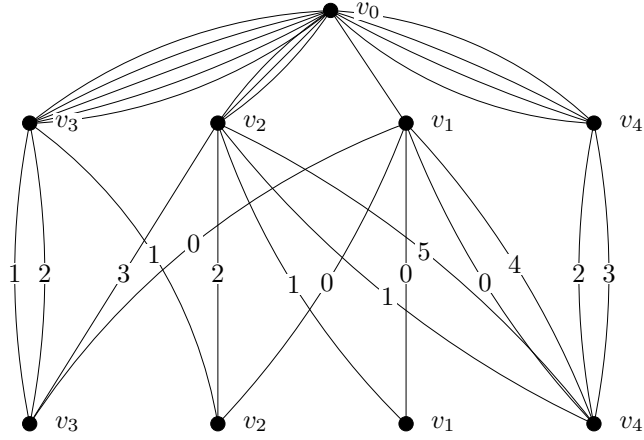
\unitlength = 0.5cm
\begin{graph}(20,14)
\graphnodesize{0.4}
 \roundnode{V0}(10,13)\nodetext{V0}(1,0){$v_0$}
\roundnode{V31}(2,10)\roundnode{V21}(7,10)\roundnode{V11}(12,10)\roundnode{V41}(17,10)
\roundnode{V32}(2,2)\roundnode{V22}(7,2)\roundnode{V12}(12,2)\roundnode{V42}(17,2)
%
\nodetext{V11}(1,0){$v_1$}\nodetext{V21}(1,0){$v_2$}
\nodetext{V31}(1,0){$v_3$}\nodetext{V41}(1,0){$v_4$}
\nodetext{V12}(1,0){$v_1$}\nodetext{V22}(1,0){$v_2$}
\nodetext{V32}(1,0){$v_3$}\nodetext{V42}(1,0){$v_4$}
%
%
\edge{V0}{V11} \bow{V0}{V21}{0.05}\bow{V0}{V21}{-0.05}
\edge{V0}{V21}
\bow{V0}{V41}{0.025}\bow{V0}{V41}{0.09}
\bow{V0}{V41}{-0.025}\bow{V0}{V41}{-0.09}
 \edge{V0}{V31}
 \bow{V0}{V31}{0.04}\bow{V0}{V31}{0.08}
\bow{V0}{V31}{-0.04}\bow{V0}{V31}{-0.08}
%
\bow{V12}{V21}{0.06} \bowtext{V12}{V21}{0.06}{1}
  \edge{V12}{V11}\edgetext{V12}{V11}{0}
%
\bow{V22}{V11}{-0.05}\bowtext{V22}{V11}{-0.05}{0}

\edge{V22}{V21} \edgetext{V22}{V21}{2}

\bow{V22}{V31}{-0.1} \bowtext{V22}{V31}{-0.1}{1}
%
\bow{V32}{V11}{0.08}\bowtext{V32}{V11}{0.08}{0}

\edge{V32}{V21}\edgetext{V32}{V21}{3}

\bow{V32}{V31}{0.05}\bow{V32}{V31}{-0.05}
\bowtext{V32}{V31}{0.05}{1}\bowtext{V32}{V31}{-0.05}{2}
%
\bow{V42}{V41}{0.05}\bow{V42}{V41}{-0.05}
\bowtext{V42}{V41}{0.05}{2}\bowtext{V42}{V41}{-0.05}{3}

\bow{V42}{V11}{0.05}\bow{V42}{V11}{-0.05}
\bowtext{V42}{V11}{0.05}{0}\bowtext{V42}{V11}{-0.05}{4}

\bow{V42}{V21}{0.06}\bow{V42}{V21}{-0.06}
\bowtext{V42}{V21}{0.06}{1}\bowtext{V42}{V21}{-0.06}{5}

\end{graph}

\medbreak
%
%

\section*{Appendix: Description of the phase space  $X_\sigma$}

%
%
%
\setcounter{theorem}{0}
\renewcommand{\thetheorem}{A.\arabic{theorem}}
\setcounter{equation}{0}
\renewcommand{\theequation}{A.\arabic{equation}}

Here we give  a combinatorial description of the phase space
$X_\sigma$ of a substitutional dynamical system assuming that
$|\sigma^n(a)|\to\infty$ for all $a\in A$.

\begin{itemize}

\item[A.] Denote by $\Lambda$ the set of all sequences $\overline
s=\{(a_n,b_n)\}_{n\geq 0}$ such that $a_nb_n\in L(\sigma)$ and
$a_n$ appears at the $(|\sigma(a_{n+1})|-1)$-th position of
$\sigma(a_{n+1})$ and $b_n$ appears at zero position of
$\sigma(b_{n+1})$.

\item[B.] Denote by $\mathcal M$ the set of all sequences
$\overline m=\{(a_0,i_0),(a_1,i_1),\ldots\}$ such that $a_j\in A$
and $i_j$ is a place of occurrence of $a_{j-1}$ in $\sigma(a_j)$,
$i_j\in \{0,1,\ldots,|\sigma^n(a_j)|-1\}$. We assume that $i_0=0$.

For a fixed $\overline m$, define inductively a sequence $\{j_n\}$
as follows: $j_0=i_0=0$ and

\begin{equation}\label{Definition_CuttingPoint}j_{n+1}=\left\{\begin{array}{ll}
|\sigma^n(a_0^{(n+1)})|+\ldots +
|\sigma^n(a^{(n+1)}_{i_{n+1}-1})|+j_n & \mbox{if }i_{n+1}\geq 1, \\
j_n & \mbox{if } i_{n+1}=0,
\end{array}\right.
\end{equation}
here $\sigma(a_{n+1})=a_0^{(n+1)}a_1^{(n+1)}\ldots
a^{(n+1)}_{|\sigma(a_{n+1})|-1}$ with $a^{(n+1)}_i\in A$.

\item[C.] Denote by $\mathcal M_0$ the set of all $\overline m$
for which $j_n\to\infty$ and $(|\sigma^n(a_n)|-j_n)\to\infty$ as
$n\to\infty$.
\end{itemize}

\begin{center}\it Construction of sequences from $X_\sigma$
\end{center}

\begin{itemize}

\item For every sequence $\overline s=\{(a_n,b_n)\}_{\geq
0}\in\Lambda$, define the sequence $w=w(\overline s)\in X_\sigma$ by
\begin{equation}
w_{[-|\sigma^n(a_n)|,|\sigma^n(b_n)|-1]}=\sigma^n(a_n).\sigma^n(b_n)\mbox{
for all }n\geq 0.
\end{equation}

\item For every $\overline m\in\mathcal M_0$, define the sequence
$w=w(\overline m)$ by

\begin{equation}
w_{[-j_n,|\sigma^n(a_n)|-j_n-1]}=\sigma^n(a_n)
\end{equation}
\end{itemize}
(it is easily seen that the sequences $w(\overline s)$ and
$w(\overline m)$ indeed belong to $X_\sigma$).

\medskip Let us recall that $T_\sigma : A^\mathbb Z\rightarrow A^\mathbb Z$ denote the shift. The next theorem describes the phase space of
an arbitrary substitution $\sigma$.

%
%

\begin{theorem} Let $\sigma: A\rightarrow A^+$ be a substitution with
$|\sigma^n|\to\infty$.  (1)
 \begin{equation}\label{Equation_RepresentationSequences}X_\sigma=\{w(\overline m) : \overline m\in\mathcal M_0\}\cup
 \bigcup_{\overline s\in\Lambda}Orb_{T_\sigma}(w(\overline s))\end{equation}

 (2) the set  $\{\omega(\overline{s}):\overline{s} \in \Lambda\}$ is finite.

 (3) If $\sigma$ is aperiodic, then $$\bigcap_{n=0}^\infty
\sigma^n(X_\sigma)=\{w(\overline s)\;:\; \overline s\in
\Lambda\}.$$

\end{theorem}
\proc{Proof.} (1) It  follows from Remark \ref{Remark_1-cutting} that
  for every $t\in X_\sigma$, there exists a sequence
$\mathcal F=\{(x_n,(y_n,i_n),z_n)\}_{n\geq 0}$ such that
\begin{itemize}
\item[i.] $x_n,y_n,z_n\in A$ and $x_ny_nz_n\in L(\sigma)$;

\item[ii.]  $i_n\in \{0,\ldots,|\sigma^n(y_n)|-1\}$ is the position
of occurrence of $y_{n-1}$ within $\sigma(y_n)$;

\item[iii.] If $i_n=0$, then $x_{n-1}$ appears in $\sigma(x_n)$ at
the position $|\sigma(x_n)|-1$;

\item[iv.] If $i_n=|\sigma(y_n)|-1$, then $z_{n-1}$ appears in
$\sigma(z_n)$ at zero position;

\item[v.] If $0< i_n<|\sigma(y_n)|-1$, then
$x_{n-1}y_{n-1}z_{n-1}$ appears in $\sigma(y_n)$ at the position
$i_n-1$.

\item[vi.] If the sequence $\{j_n\}$ is determined by the sequence
$\overline m =\{(y_n,i_n)\}_{n\geq 0}\in\mathcal M$ as in
(\ref{Definition_CuttingPoint}), then
\begin{equation}\label{Claim2_Equality}
t{[-j_n-|\sigma^n(x_n)|,|\sigma^n(y_n)|-j_n-1+|\sigma^n(z_n)|]}=
\sigma^n(x_n)\sigma^n(y_n)\sigma^n(z_n).
\end{equation}
\end{itemize}

Now we have three options:

(a) $j_n\to\infty$ and $|\sigma^n(y_n)|-j_n\to\infty$ as
$n\to\infty$. In this case, we get that $t=w(\overline m)$.

(b) The sequence $\{j_n\}$ is bounded and
$(|\sigma^n(y_n)|-j_n)\to\infty$. Here, we get that $w=w(\overline m)$ is one-sided (to
the right). Since the sequence $\{j_n\}$ is non-decreasing and
bounded, there is $n_0>0$ such that $j_n=j_{n_0}$ for all $n\geq
n_0$. This implies that $i_n=0$ for all $n\geq n_0$, i.e. $y_n$
appears at zero position in $\sigma(y_{n+1})$. This, in
particular, implies that $x_n$ appears in $\sigma(x_{n+1})$ at the
last position. Therefore,
$$t{[-j_{n_0}-|\sigma^n(x_n)|,|\sigma^n(y_{n})|-j_{n_0}-1]}=
\sigma^n(x_n)\sigma^{n}(y_{n})
\mbox{ for }n\geq n_0.$$ Set $x_{n_0-1}'=\sigma(x_{n_0})[|\sigma(x_{n_0})|-1]$ and $y_{n_0-1}'=\sigma(y_{n_0})[0]$. Define inductively $x_{k-1}'=\sigma(x_{k}')[|\sigma(x_k')|-1]$ and $y_{k-1}'=\sigma(y_{k}')[0]$ for $k=n_0-1,\ldots, 1$. Setting   $$\overline
s=\{(x_0',y_0'),\ldots,(x_{n_0-1}',y_{n_0-1}'),(x_{n_0},y_{n_0}),(x_{n_0+1},y_{n_0+1}),\ldots\},$$ we get that $t=T^{j_0}(w(\overline
s))$.

(c) The sequence $\{|\sigma^n(y_n)|-j_n\}$ is bounded and
$j_n\to\infty$ as $n\to\infty$. The proof in this case is similar to (b).

(2) Denote by $A_0$ (by $A_1$) the set of all letters $a\in A$ such that $\sigma^{n_a}(a)$ begins (ends) with $a$ for some $n_a$. Clearly, if such an $n_a$ exists, then it can be chosen from the interval $[1,|A|+1]$. Set $p=\prod_{a\in A_0}n_a\times \prod_{b\in A_1}n_b$. Then $p$ is bounded by $(|A|+1)^{2|A|}$.

Consider a sequence $\overline s=\{(a_n,b_n)\}_{n\geq 0}\in\Lambda$. Find an infinite set $I$ such that $a_k=a_{k'}$ and $b_k=b_{k'}$ for all $k,k'\in I$. Therefore, $\sigma^{k-k'}(a_k)$ ends with $a_k$ for all $k>k'$, $k,k'\in I$. This shows that $\sigma^p(a_k)$ ($\sigma^p(b_k)$) ends (begins) with $a_k$ ($b_k$) for all $k\in I$.
We can write
down each $k\in I$ as
$$k=pl_k+m_k\mbox{ with }l_k\geq 0\mbox{ and }0\leq
m_k<p.$$ Find an infinite set $J\subseteq I$ such that
$m_k=m_{k'}$ for $k,{k'}\in J$. It follows that $a_{m_k+pn}=a_{m_k}$ and $b_{m_k+pn}=b_{m_k}$ for all $n\geq 0$ and $k\in J$.
Then $$w(\overline s)=\lim\limits_{n\to\infty}\sigma^n(a_n).\sigma^n(b_n)=
\lim\limits_{n\to\infty}\sigma^{pn}
(\sigma^{m_k}(a_{m_k})).\sigma^{pn}(\sigma^{m_k}(b_{m_k})),
$$ where $k$ is any integer from $J$. This shows that each element $w(\overline s)$ is determined by a finite number of parameters taken from finite sets.

 (3) If $y\in Y=\bigcap_{n\geq 0}\sigma^n(X_\sigma)$, then
for every $n\geq 0$ we have $y=\sigma^n(x_n)$ for some $x_n\in
X_\sigma$. By Theorem \ref{TheoremRecognizability}, the point
$x_n$ is   uniquely defined. Therefore, $\sigma(x_{n+1})=x_n$
for every $n$. Setting $\overline s=\{(x_{n}[-1],x_n[0])\}_{n\geq
0}$, we get that $y=w(\overline s)$.

Conversely, if $y=w(\overline s)$, then it is not hard to
decompose $y$ into $n$-words such that $y=\sigma^n(x_n)$ for some
$x_n\in X_\sigma$. In particular, this shows that $y\in Y$.\ep
\medskip

{\it Acknowledgments}. The work was done during our mutual visits to the University of Torun and Institute for Low Temperature Physics. We are thankful to these institutions for the hospitality and support. We would like to express our special thanks to  C.~Skau and B.~Solomyak  for numerous stimulating discussions and to T. Downarowicz and A. Maass for sending us their paper and drawing our attention to their result. K.~Medynets also acknowledges  the support  of INTAS  Young Scientist Fellowship.



%


\end{document}